\documentclass[11pt]{article}

\PassOptionsToPackage{numbers, compress}{natbib}
\usepackage[nonatbib,preprint]{neurips_2019}

\usepackage{amssymb}
\usepackage{amsmath}
\usepackage{amsthm}
\usepackage[utf8]{inputenc} % allow utf-8 input
\usepackage[T1]{fontenc}    % use 8-bit T1 fonts
\usepackage{url}            % simple URL typesetting
\usepackage{booktabs}       % professional-quality tables
\usepackage{amsfonts}       % blackboard math symbols
\usepackage{nicefrac}       % compact symbols for 1/2, etc.
\usepackage{microtype}      % microtypography
\usepackage{graphicx}
\usepackage{epstopdf}
\usepackage{subcaption}
\usepackage{hyperref}
\usepackage{xcolor}
\usepackage{algorithm,algpseudocode}
\usepackage{natbib}

\newtheorem{theorem}{Theorem}
\newtheorem{rem}{Remark}

\title{Fast approximation of orthogonal matrices and application to PCA}
\title{Fast approximation of orthogonal matrices\\ and application to PCA}

\author{%
  Cristian Rusu%\thanks{Use footnote for providing further information
\\
Faculty of Automatic Control and Computers, University Politehnica Bucharest\\
  \texttt{cristian.rusu@upb.ro} \\
   Lorenzo Rosasco \\
LCSL, Universit\'a di Genova\\
 Massachusetts Institute of Technology and Istituto Italiano di Tecnologia
}

\begin{document}

\maketitle

\begin{abstract}
	We study the problem of approximating orthogonal matrices  so that their application is numerically  fast and yet accurate. We find an approximation by solving an optimization problem over a set of structured matrices, that we call extended orthogonal Givens transformations, including Givens rotations as a special case. We propose an efficient greedy algorithm to solve such a problem and show that it strikes a balance between approximation accuracy and speed of computation. The approach is relevant to spectral methods and we illustrate its application to PCA.
\end{abstract}
\section{Introduction}

Orthonormal transformations play a key role in most matrix decomposition techniques and spectral methods \cite{Belabbas369}. As such, manipulating them in an efficient manner is essential to many practical applications. While general matrix-vector multiplications with orthogonal matrices take $O(d^2)$ space and time, it is natural to ask whether faster approximate computations (say $O(d \log d)$) can be achieved while retaining enough accuracy.\\
Approximating an orthonormal matrix with just a few building blocks is hard in general. The standard decomposition technique meant to reduce complexity is a low-rank approximation. Unfortunately, for an orthonormal matrix, which is perfectly conditioned, this approach is meaningless.\\
In this work, we are inspired by the fact that several orthonormal/unitary transformations that exhibit low numerical complexity are known. The typical example is the discrete Fourier transform with its efficient implementation as the fast Fourier transform \cite{doi:10.1137/1.9781611970999} together with other Fourier-related algorithms: fast Walsh-Hadamard transforms \cite{1674569}, fast cosine transforms \cite{1163351}, and fast Hartley transforms \cite{1457236}. Other approaches include fast wavelet transforms \cite{doi:10.1002/cpa.3160440202}, banded orthonormal matrices \cite{SIMON2007120, Strang12413} and fast Slepian transforms \cite{KARNIK2019624}. Decomposition of orthogonal matrices into $O(d)$ Householder reflectors or $O(d^2)$ Givens rotations \cite{Golub1996}[Chapter 5.1] are known already. These basic building blocks have been further extended, for example we have fast Givens rotations \cite{Rath1982} and two generalizations of the Givens rotations \cite{BILOTI201356} and \cite{GGR}. In theoretical physics, unitary decompositions parametrize symmetry groups \cite{Tilma:2002ke}, and they are compactly parametrized using $\sigma$-matrices \cite{Spengler_2010} or symmetric positive definite matrices \cite{Barvinok}. To the best of our knowledge, none of these factorizations focus on reducing the computational complexity of using the orthogonal/unitary transformations but rather they model properties of physical systems.\\
Our idea is to approximately factor any orthonormal matrix into a product of a fixed number of sparse matrices such that their application (to a vector) 
%computational  of matrix-vector multiplication 
has linearithmic complexity. In this paper, we derive structured approximations to orthonormal matrices that can be found efficiently and applied remarkably fast. We pose the search for an efficient approximation as an optimization problem over a product of structured matrices, the extended orthogonal Givens transformations. These structures extend Givens rotations to also include reflectors with no computational drawback and suggest a decomposition of the main optimization problem into sub-problems that are easy to understand and solved via a greedy approach. The theoretical properties of the obtained solution are characterized in terms of approximation bounds while the empirical properties are studied extensively.\\
We illustrate our approach  considering dimensionality reduction with principal component analysis (PCA). %We show the application of our proposed method to the calculation of fast principal component analysis (PCA) projections.\\
Here the goal is not to propose a new fast algorithm to compute the principal directions, plenty of efficient algorithms  for this task \cite{FastMCSVD}\cite{Golub1996}[Chapter~10] \cite{FindingStructureWithRandomness}. Rather, we aim at constructing fast dimensionality reduction operators.
%PCA provides a fundamental solution to the problem of dimensionality reduction with projections that are excellent candidates for the fast approximations we consider. In most machine learning applications, using PCA requires two steps. First, the data matrix is factorized by the singular value decomposition with the Lanczos method \cite{Golub1996}[Chapter~10] or, randomized methods \cite{FastMCSVD, FindingStructureWithRandomness}. Second, the singular vectors (principal components) define a projection to reduce data dimensionality. In this paper, we focus on this second step.\\
While the calculation of the principal components is a one-off computation, a numerically efficient projection operation is critical since it is required multiple times in downstream applications. The problem of deriving fast projections has also been previously studied. Possible approaches include: fast wavelet transforms \cite{WaveletsForPCA}, sparse PCA \cite{MINIMAXOPTIMALSPCA, SparsePCA}, structured transformations such as circulant matrices \cite{CirculantDRoperator}, Kronecker products \cite{KroneckerPCA}, Givens rotations \cite{DCTandKLT, Treelets, LeCunFastApproximations} or structured random projections \cite{FastJL, ToeplitzJL}. Compared to these works, we propose a new way to factorize any orthogonal matrix, including PCA directions, into simple orthogonal structures that we call extended orthogonal Givens transformations and which naturally lead to optimization problems that have closed-form solutions and are therefore efficiently computed.\\
We note that our approach provides new perspectives on the structure of the orthogonal group and how to coarsely approximate it, which might have an impact on other open research questions.\\% where orthogonal matrices are involved.\\
% The practical purpose of our work is to close the computational gap between the optimal, but unstructured, PCA projections and the classic fast, but not data dependent, transformations enumerated above. This would extend the application of PCA projections to scenarios where we are dealing with high-dimensional data or where the running time (numerical complexity) and energy consumption are factors of critical importance.\\
The paper is organized as follows: Section 2 describes the basic building blocks and algorithm that we propose, Section 3 gives the theoretical guarantees for our contributions, Section 4 details the application of our method to PCA projections and Section 5 shows the numerical experiments.

\section{The proposed algorithm}
	Given a $d\times d$ orthonormal $\mathbf{U}$, the matrix-vector multiplication $\mathbf{Ux}$ takes $O(d^2)$ operations. We want to build $\mathbf{\bar{U}}$ such that $\mathbf{U} \approx \mathbf{\bar{U}}$ and $\mathbf{\bar{U}x}$ takes $O(d \log d)$ operations. Parametrizations of orthonormal matrices \cite{doi:10.1137/0908055} are known, but to be best of our knowledge, the problem of accurately approximating $\mathbf U$ as  product of only a few $O(d \log d)$ transformations is open. Given a $d \times p$ diagonal $\mathbf{\Sigma}_p$, in the spirit of previous work minimizing Frobenius norm approximations \cite{Kondor2014MMF, lemagoarou:hal-01416110}, we consider the problem 
	\begin{equation}
	\underset{\mathbf{\bar{U}}, \ \mathbf{\bar{\Sigma}}_p}{\text{minimize}}\ \| \mathbf{U} \mathbf{\Sigma}_p - \mathbf{\bar{U}}\mathbf{\bar{\Sigma}}_p  \|_F^2 \text{ subject to } \mathbf{\bar{U}} \in \mathcal{F}_g,
	\label{eq:theoptimizationproblem_only_orthonormal}
	\end{equation}
	where $\mathbf{\bar{U}}$ is $d \times d$ and $\mathbf{\bar{\Sigma}}_p$ is a $d \times p$ diagonal matrix. Choosing $p = d$ while $\mathbf{\Sigma}_p$ and $\mathbf{\bar{\Sigma}}_p$ to be the identity, we simply approximate $\mathbf U$. We will also use $\mathbf{U}_p$ do denote the first $p \leq d$ columns of $\mathbf{U}$.
The above general formulation allows to also consider cases where  different directions might have different importance. 
Then, $\mathcal{F}_g$ is a  set of orthogonal matrices -- defined next -- that can be applied fast and allow to efficiently but approximately solve
%matrix-vector multiplication is of order $O(g)$ and 
\eqref{eq:theoptimizationproblem_only_orthonormal}.
% can be  approximately solved.

\subsection{The basic building blocks}

Classic matrix building blocks that are numerically efficiency include  circulant/Toeplitz matrices or Kronecker products. These choices are inefficient as they depend on $O(d)$ free parameters but their matrix-vector product cost is $O(d \log d)$ or even $O(d \sqrt{d})$, i.e., they do not scale linearly with the number of parameters they have. Consider the sparse orthogonal matrices
\begin{equation}
\mathbf{G}_{ij} = \begin{bmatrix} \mathbf{I}_{i-1} &  & &\\
& * & & *  \\
& & \mathbf{I}_{j-i-1} & \\
& * &  & *  \\
& & & & \mathbf{I}_{d-j} \\
\end{bmatrix}, \mathbf{\tilde{G}}_{i j} \in \left\{ \! \begin{bmatrix} c & -s \\ s & c \end{bmatrix}\!,\! \begin{bmatrix} c & s \\ s & -c	\end{bmatrix}   \! \right\},\text{ such that } c^2 + s^2 = 1,
\label{eq:theG}
\end{equation}
where the non-zero part (denoted by $*$ and $\mathbf{\tilde{G}}_{i j}$) on rows and columns $i$ and $j$. The transformation in \eqref{eq:theG}, with the first option in $\mathbf{\tilde{G}}_{i j}$, is a Givens (or Jacobi) rotation. With the second option, we have a very sparse Householder reflector. These transformations were first used by \citet{FastSparsifyingTransforms} to learn numerically efficient sparsifying dictionaries for sparse coding. The $\mathbf{G}_{ij}$s have the following advantages: i) they are orthogonal; ii) they are sparse and therefore fast to manipulate: matrix-vector multiplications $\mathbf{G}_{ij}\mathbf{x}$ take only $6$ operations; iii) there are two degrees of freedom to learn: $c$ (or $s$) and the binary choice; and iv) allowing both sub-matrices in $\mathbf{\tilde{G}}_{i j}$ enriches the structure and as we will see, leads to an easier (closed-form solutions) optimization problem.\\
We propose to consider matrices $\mathbf{\bar{U}} \in \mathcal{F}_g$ that are   products of $g$ transformations from \eqref{eq:theG}, that is
\begin{equation}
\mathbf{\bar{U}} =
 \prod_{k=1}^{g} \mathbf{G}_{i_k j_k} = \mathbf{G}_{i_1 j_1}  \dots \mathbf{G}_{i_g j_g}.
\label{eq:approxu}
\end{equation}
Matrix-vector multiplication with $\mathbf{\bar{U}}$ takes $6g$ operations -- when $g$ is $O(d \log d)$ this is significantly better than $O(d^2)$, while the coding complexity of each $\mathbf{G}_{i j}$ is approximately $2\log_2 d + C$: $2\log_2 d -1$ bits to encode the choice of the two indices, a constant factor $C$ for the pair $(c,s)$ and 1 bit for the choice between the rotation and reflector. The coding complexity of $\mathbf{\bar{U}}$ scales linearly with $g$.\\
	We note that Givens rotations have been used extensively to build numerically efficient transformations \cite{5560826, EffectiveGivens, Kondor2014MMF, lemagoarou:hal-01416110, Treelets}. However,  $2 \times 2$ reflector was not used before. This may be because in linear algebra (e.g. in QR factorization) and in optimization \cite{Shalit14} considering also the reflector has no additional benefit: the rotation alone introduces each zero in the QR factorization and the reflector does not have an exponential mapping on the orthogonal manifold, respectively. As we will show, considering both the rotation and the reflector has the advantage of providing a closed-form solution to our problem.

\subsection{The proposed greedy algorithm}

We propose to solve the optimization problem in \eqref{eq:theoptimizationproblem_only_orthonormal} with a greedy approach: we keep $\mathbf{\bar{\Sigma}}_p$ and all variables fixed except for a single $\mathbf{G}_{i_k j_k}$ from $\mathbf{\bar{U}}$ and minimize the objective function. When optimizing w.r.t.  $\mathbf{G}_{i_k j_k}$ it is convenient to write 
% \prod_{k=1}^{g} \mathbf{G}_{i_k j_k} 
\begin{equation}
\! \| \mathbf{U}\mathbf{\Sigma}_p  \! \! - \! \mathbf{\bar{U} \bar{\Sigma}}_p \|_F^2 %\| \mathbf{U}_p\mathbf{\Sigma}_p - \mathbf{G}_{i_1 j_1}  \dots \mathbf{G}_{i_g j_g} \mathbf{\bar{\Sigma}}_p \|_F^2  %\\
\!  = \!  \| 
%\mathbf{G}_{i_{k-1} j_{k-1}}^T \dots \mathbf{G}_{i_1 j_1}^T  
\underbrace{ \prod_{t=1}^{k-1} \mathbf{G}_{i_t j_t}^T   
\mathbf{U} \mathbf{\Sigma}_p \! \!  }_{{\mathbf L}^{(k)}}
- \mathbf{G}_{i_k j_k} 
\underbrace{
 \prod_{t=k+1}^{g} \mathbf{G}_{i_t j_t}
%\dots \mathbf{G}_{i_g j_g}
 \mathbf{\bar{\Sigma}}_p}_{{\mathbf N}^{(k)}} \|_F^2 \!   = \! \| \mathbf{L}^{(k)} \! - \! \mathbf{G}_{i_k j_k} \mathbf{N}^{(k)} \|_F^2.
\label{eq:thestart}
\end{equation}
The next result characterizes the Givens transformation $\mathbf{G}_{i_k j_k}$ minimizing the above norm. We drop the dependence on $k$ for ease of notation.

%\begin{equation}
%\! \| \mathbf{U}_p\mathbf{\Sigma}_p  \! \! - \! \mathbf{\bar{U} \bar{\Sigma}}_p \|_F^2 =  %\| \mathbf{U}_p\mathbf{\Sigma}_p - \mathbf{G}_{i_1 j_1}  \dots \mathbf{G}_{i_g j_g} \mathbf{\bar{\Sigma}}_p \|_F^2  %\\
%\! \! = \! \!  \| \mathbf{G}_{i_{k-1} j_{k-1}}^T \dots \mathbf{G}_{i_1 j_1}^T  \mathbf{U}_p \mathbf{\Sigma}_p \! \!  - \! \! \mathbf{G}_{i_k j_k} \dots \mathbf{G}_{i_g j_g} \mathbf{\bar{\Sigma}}_p \|_F^2 \! \!  = \!\! \| \mathbf{L}^{(k)} \! - \! \mathbf{G}_{i_k j_k} \mathbf{N}^{(k)} \|_F^2.
%\label{eq:thestart}
%\end{equation}

\begin{theorem}
%\noindent \textbf{Theorem 1 
[Locally optimal $\mathbf{G}_{ij}$] 
Let $\mathbf{L}$ and $\mathbf{N}$ be two $d \times p$ matrices. Further, let  $\mathbf{Z} = \mathbf{L} \mathbf{N}^T$ and $\mathbf{Z}_{\{i,j\}} = \begin{bmatrix} Z_{ii} & Z_{ij}  \\ Z_{ji} & Z_{jj} \end{bmatrix}$ then we have
\begin{equation}
\ C_{ij} = \| \mathbf{Z}_{ \{i,j\} } \|_* - \text{tr}(\mathbf{Z}_{ \{ i,j \} }) = \left\{ \begin{array}{lr} \sqrt{ (Z_{ii} + Z_{jj})^2 + (Z_{ij} - Z_{ji})^2} -Z_{ii} -Z_{jj} , \text{if } \det( \mathbf{Z}_{ \{ i,j \} } )  \geq 0 \\ \sqrt{ (Z_{ii} - Z_{jj})^2 + (Z_{ij} + Z_{ji})^2} -Z_{ii} -Z_{jj}, \text{if } \det( \mathbf{Z}_{ \{ i,j \} } )  < 0 \end{array}\right.
\label{eq:theCijcomputed}
\end{equation}
Let $\mathbf{Z}_{\{i^\star,j^\star\}} \! = \! \mathbf{V}_1 \mathbf{S  V}
_2^T$ be the SVD of $\mathbf{Z}_{\{i^\star,j^\star\}}$ , with
\begin{equation}
(i^\star, j^\star) = \underset{(i, j),\ j > i}{\arg \max}\ \ C_{ij}.
\label{eq:theCij}
\end{equation}
Then, the  extended orthogonal Givens transformation that minimizes $\| \mathbf{L} - \mathbf{G}_{i j} \mathbf{N} \|_F^2$ is given by
$\mathbf{\tilde{G}}_{i^\star j^\star}^\star = \mathbf{V}_1 \mathbf{V}_2^T.$
\end{theorem}

%  where $\mathbf{Z} = \mathbf{L} \mathbf{N}^T$ and the optimal indices
%\begin{equation}
%(i^\star, j^\star) = \underset{(i, j),\ j > i}{\arg \max}\ \ C_{ij}, \ C_{ij} = \| \mathbf{Z}_{ \{i,j\} } \|_* - \text{tr}(\mathbf{Z}_{ \{ i,j \} })\text{ where } \mathbf{Z}_{\{i,j\}} = \begin{bmatrix} Z_{ii} & Z_{ij}  \\ Z_{ji} & Z_{jj} \end{bmatrix}.
%\label{eq:theCij}
%\end{equation}
The above theorem  derives a locally optimal way to construct an approximation $\mathbf{\bar{U}}$. We iteratively apply the result to find, for each component $k$ in \eqref{eq:approxu}, the extended orthogonal Givens transformation that best minimizes the objective function \eqref{eq:theoptimizationproblem_only_orthonormal}. The full procedure is  in Algorithm 1 and can be viewed in two different ways: i) a coordinate minimization algorithm; or ii) a hierarchical decomposition where each stage is extremely sparse. The proposed algorithm is guaranteed to converge, in the sense of the objective function \eqref{eq:theoptimizationproblem_only_orthonormal}, to a stationary point. Indeed, no step in the algorithm can increase the objective function, since the sub-problems are minimized exactly: we choose the best indices and then perform the best $2 \times 2$ transformation. We note three remarks on the properties of Algorithm 1.
\begin{rem}[Complexity of Algorithm 1]
The computational complexity of the iterative part of Algorithm 1 is $O(d g)$ and the initialization is dominated by the computation of all the scores $C_{ij}$ which takes $O(d^2)$. Note that, the $C_{ij}$s are computed from scratch only once in the initialization phase. After that,  at each step $k$ we need to  recompute the $C_{ij}$ (redo the $2\times 2$ singular value decompositions) only for the indices $(i_k, j_k)$  currently used (the Givens transformations act on two coordinates at a time). All other scores are update by the same quantity: in \eqref{eq:theCij}, the $C_{ij}$ are the same except when $i$ or $j$ belong to the set $(i_k, j_k)$. This observation substantially reduces the running time.
\end{rem}
\begin{rem}[Complexity of applying $\mathbf{\bar{U} \mathbf{\bar{\Sigma}}}_p$] 
When  $p <   d$ the computational complexity of $6g$ operations is an upper bound. Since we keep only  $p$ components, we need  be careful  not to perform operations whose result is thrown away by the mask $\mathbf{\bar{\Sigma}}_p$. Consider for example a transformation $\mathbf{G}_{1d}$ applied to a vector of size $d$  projected to a $p < d$ dimensional space. The three operations that take place on the $d^\text{th}$ component are unnecessary. Then, after  computing $\mathbf{\bar{U}}$, a pass is made through each of the $g$ transforms to decide which of two coordinates the computations are necessary for the final result. As we will show, this further improves the numerical efficiency of our method.
\end{rem}
\begin{rem}[On the choice of indices] Algorithm 1 greedily chooses at each step $k$ the indices  according to \eqref{eq:theCij}. Other factors might be considered: i) choosing indices based on previous choices so that only a select group of indices are used throughout the algorithm, or ii) make multiple choices at each step in order to speed up the algorithm.
\end{rem}
\begin{algorithm}[tb]
	\caption{Approximate orthonormal matrix factorization with extended Givens transformations}
	\label{alg:fpca}
	\begin{algorithmic}
		\State {\bfseries Input:} The $p$ orthogonal components $\mathbf{U}_p$ and their weights (singular values) $\mathbf{\Sigma}_p$, the size $g$ of the approximation \eqref{eq:approxu}, the update rule for $\mathbf{\bar{\Sigma}}_p$ in \{ `identity', `original', `update' \} and the stopping criterion $\epsilon$ (default taken to be $\epsilon = 10^{-2}$).
		
		\State {\bfseries Output:} The linear transformation $\mathbf{\bar{U}} \mathbf{\mathbf{\bar{\Sigma}}}_p$, the approximate solution to \eqref{eq:theoptimizationproblem_only_orthonormal}.
		\State {\bfseries Initialize:} $\mathbf{G}_{i_k j_k} \! = \! \mathbf{I}_{d \times d},k\! =\! 1,\dots,g$ and compute all scores $C_{ij}$ according to \eqref{eq:theCijcomputed} with $\mathbf{Z} = \mathbf{U}_p \mathbf{\bar{\Sigma}}_p^T$, where $\mathbf{\bar{\Sigma}}_p = \begin{bmatrix} \mathbf{I}_{p \times p}; & \mathbf{0}_{(d-p) \times p} \end{bmatrix}$ if the update rule is `identity' and $\mathbf{\bar{\Sigma}}_p = \mathbf{\Sigma}_p$ otherwise.
		\Repeat
		\State Set $\mathbf{L}^{(0)} = \mathbf{U}_p \mathbf{\Sigma}_p$ and set $\mathbf{N}^{(0)} = \mathbf{G}_{i_1 j_1} \dots \mathbf{G}_{i_g j_g} \mathbf{\bar{\Sigma}}_p$.
		
		\For{$k = 1$ {\bfseries to} $g$ }
		\State Update $\mathbf{N}^{(k)} = \mathbf{G}_{i_k j_k}^T \mathbf{N}^{(k-1)}$ and find best score according to \eqref{eq:theCij}.
		\State Compute the best $k^\text{th}$ transformation by Theorem 1.
		\State Update $\mathbf{L}^{(k)} = \mathbf{G}_{i_k j_k}^T \mathbf{L}^{(k-1)}$ and then update all scores $C_{ij}$ in \eqref{eq:theCijcomputed} but only for indices $i, j \in \{i_k, j_k\}$ with $\mathbf{Z}\! = \! \mathbf{L}^{(k)} (\mathbf{N}^{(k)})^T$ -- all other scores are unchanged.
		\EndFor
		
		\State Set $\mathbf{\bar{\Sigma}}_p = \begin{bmatrix}
		\text{diag}(\mathbf{L}^{(g)}); & \mathbf{0}_{(d-p) \times p}
		\end{bmatrix}$ if rule is `update', $i \leftarrow i+1$ and $\epsilon_i = \| \mathbf{L}^{(g)} - \mathbf{\bar{\Sigma}}_p \|_F^2 $.
		\Until $ |\epsilon_{i-1} - \epsilon_i| < \epsilon$, if $i > 1$.
	\end{algorithmic}
\end{algorithm}

\section{Analysis of the proposed algorithm}

%In this section we explore bounds on the performance of the proposed algorithm.
We consider $p = d$, i.e., $\mathbf{\bar{\Sigma}}_p = \mathbf{I}_{d \times d}$ and therefore $\mathbf{Z} = \mathbf{U}$. We model the $\mathbf{U}$ as a random orthonormal matrix with Haar measure \cite{HowToRandomUnitary} updated so that the diagonal is positive. We perform this update because multiplication by a diagonal matrix with $\pm 1$ entries has no computational cost but it brings $\mathbf{U}$ closer to $\mathbf{I}_{d \times d}$.
The goal of this section is to establish upper bounds for the distance between $\mathbf{U}$ and $\mathbf{\bar{U}}$, as a function of $d$ and $g$. We first comment  on the inherent difficulty of the problem.

\begin{rem}[The approximation gap] Since  the orthogonal group has size $O(d^2)$, by the pigeonhole principle  a random orthogonal matrix as \eqref{eq:approxu} and only $g \ll d^2$ degrees of freedom cannot be exactly  approximated with less than $O(d^2)$ operations. 
For our purposes, think $g$ either $O(d)$ or $O(d \log d)$. 
%This suggests  we cannot exactly approximate any orthonormal matrix %It is therefore unreasonable to try to improve the complexity of matrix-vector multiplication in general, while the matrix-matrix multiplication problem is still open \cite{LeGall:2014:PTF:2608628.2608664}. 
Our goal is to show that the fast structures we propose can perform well in practice and have theoretical bounds that guarantee worse case or average accuracy.
\end{rem}
\noindent Next, we show two approximations bounds depending on the number of Givens transformations \eqref{eq:theG}.
\begin{theorem}[A special bound] Given a random $d \times d$ orthonormal $\mathbf{U}$, for large $d$, its approximation $\mathbf{\bar{U}}$ from \eqref{eq:theoptimizationproblem_only_orthonormal} with $g = d/2$ transformations  from \eqref{eq:theG} obeys
\begin{equation}
\mathbb{E}[ \| \mathbf{U} - \mathbf{\bar{U}} \|_F^2] \leq 2d - \sqrt{2\pi d}.
\end{equation}
\end{theorem}
\begin{theorem}[A general bound] Given a random $d \times d$ orthonormal $\mathbf{U}$, for large $d$, its approximation $\mathbf{\bar{U}}$ from \eqref{eq:theoptimizationproblem_only_orthonormal} with $g \leq d(d-1)/2$ transformations from \eqref{eq:theG} is bounded by
\begin{equation}
\mathbb{E}  \left[ \| \mathbf{U} - \mathbf{\bar{U}} \|_F^2 \right] \leq 2(d - \lfloor r \rfloor) - \frac{2\sqrt{2}}{\sqrt{\pi}} \sqrt{d - \lfloor r \rfloor},\text{ where } r = d - \frac{1 + \sqrt{(2d-1)^2 - 8g}}{2}.
\label{eq:theresultsofu3}
\end{equation}
\end{theorem}
Theorem 2 shows  that, on average, the performance might degrade with increasing $d$. As  stated in Remark 4, this is not surprising since the orthonormal group is much larger than the structure we are trying to approximate it with. The next result provides a bound for  other values of $g$. In Theorem 3, taking $g = c_1 d \log d$ for some positive constant $c_1$ we have that $r \approx c_1 \log d$. This means that whenever $p \ll d$ we will roughly need $O(d)$ Givens transformations from \eqref{eq:theG} to improve the $\lfloor r \rfloor$ term. Since the proof of the theorem uses only rotations (and furthermore, in a particular order of indices $(i_k,j_k)$) we expect our algorithm to perform much better than the bound indicates as it allows for a richer structure \eqref{eq:theG} and uses greedy steps that maximally improve the accuracy at each step.\\
The previous theorems consider  the Frobenius norm. In the Jacobi iterative process for diagonalizing a symmetric matrix with Givens rotations \cite{Golub1996}[Chapter~8.4] the progress of the procedure (convergence) is measured using the off-diagonal ``norm'' $\text{off}(\mathbf{U}) = \sqrt{\sum_{t}^d \sum_{q \neq t}^d U_{tq}^2}$.\\
\begin{theorem}[Convergence in the off-diagonal norm] Given a $d \times d$ orthonormal $\mathbf{U}$ and a single Givens transformation $\mathbf{G}_{ij}$, assuming $\text{det}(\mathbf{U}_{ \{i,j\} }) \geq 0$ we have
\begin{equation}
\text{off}(\mathbf{UG}_{ij}^T )^2 \leq \text{off}(\mathbf{U})^2 \! + \frac{1}{2}((U_{ii} - U_{jj})^2- (U_{ij} - U_{ji})^2).
\end{equation}
\end{theorem}
This result shows that, unlike with the Jacobi iterations, monotonic convergence in this quantity is not guaranteed and depends on the relative differences between the diagonal and the off-diagonal entries of $\mathbf{U}_{ \{i,j\} }$. Our method convergence monotonically to a stationary point when we measure the progress in the Frobenius norm.

\begin{rem}[The effect of a single $\mathbf{G}_{ij}$] Given a $d \times d$ orthonormal $\mathbf{U}$ and a Givens transformation $\mathbf{G}_{ij}$ from \eqref{eq:theG} we have that: i) $\mathbf{U}\mathbf{G}_{i j}^T$ is closer to the identity matrix in the sense that $\mathbf{\tilde{G}}_{ij}$ makes a positive contribution to the diagonal elements, i.e.,  $\text{tr}(\mathbf{U}\mathbf{G}_{i j}^T) = \text{tr}(\mathbf{U}) + C_{i j}$; and ii) $\mathbb{E}[C_{ij}] \approx 0.6956 d^{-1/2}$ if $\mathbf{U}$ is random with Haar measure and positive diagonal for large $d$.
\end{rem}
The above remark  suggests a metric  to study the convergence of the proposed method: each Givens transformation adds the score $C_{ij}$ to the diagonal entries of the current approximation (and therefore ensures that $\mathbf{U \bar{U}}^T $ converges to the identity -- the only diagonal orthonormal matrix). By choosing the maximum $C_{ij}$ we are taking the largest step in this direction.
\begin{rem}[Evolution of $C_{i_kj_k}$ with $k$] Given a fixed $0 < u < 1$, consider the toy construction $\mathbf{U}_{ \{i,j\} } = \begin{bmatrix} u & z_2 \\ z_1 & u \end{bmatrix}$, i.e., a $2 \times 2$ sub-matrix of a $d \times d$ orthonormal matrix where diagonal elements are equal and off-diagonals are two independent truncated standard normal random variables in the interval $[-\sqrt{1-u^2}, \sqrt{1-u^2}]$ (since rows and columns of $\mathbf{U}$ are $\ell_2$ normalized). Then, for large $d$, by direct calculation we have that the  expected score $\mathbb{E}[ C_{ij}(u) ]$, i.e., $C_{ij}$ as a function of $u$, obeys
\begin{equation}
\begin{aligned}
\mathbb{E}[C_{ij}(u)] \propto (1-u)^2,
\end{aligned}
\end{equation}
i.e., the expected $C_{ij}$ decreases on average quadratically with the increase in the diagonal elements.\\
The remark is intuitive: as $k$ increases $\mathbf{\bar{U}}$ is more accurate and $\mathbf{U\bar{U}}^T$ becomes diagonally dominant, i.e., $\mathbf{U\bar{U}}^T \to \mathbf{I}_d$ as $k \to O(d^2)$, and we do expect to reach lower scores $C_{i_k j_k}$, i.e., few Givens transformations from \eqref{eq:theG} provide a rough estimation while very good approximations require $k \approx d^2$.
\end{rem}
%\noindent Proofs of the results and details on remarks are collected in the attached supplementary materials.

\subsection{Other ways to measure the approximation error}

Throughout this paper we use the Frobenius norm to measure and study the approximation error we propose. In this section we discuss this choice and explore a few alternatives.\\
In the linear algebra literature, a natural way to measure approximation error is through the operator norm. This is especially true in the randomized linear algebra field (where matrix concentration inequalities which bound the operator norm play a central role). Moreover, the review manuscript \citep{10.1561/2200000048} deals explicitly with some potential issues that might arise from using the Frobenius, instead of the operator, norm in matrix approximations: the discussion in Chapter 6 of \citep{10.1561/2200000048} entitled ``Warning: Frobenius--Norm Bounds''. The text highlights situations where the matrix to be approximated is low rank and corrupted by noise and/or scaling issues are present. While that discussion holds true, in our case we deal with a given perfectly conditioned matrix $\mathbf{U}$ and its approximation displays the same scaling ($\mathbf{\bar{U}}$ has normalized columns just as $\mathbf{U}$). Consider the following remark.

\noindent \textbf{Remark 7. [Simplifying the Frobenius norm objective function]} \textit{Consider for simplicity $p = d$ and $\mathbf{\bar{\Sigma}}_d = \mathbf{\Sigma}_d$ in \eqref{eq:theoptimizationproblem_only_orthonormal} and see that $\| (\mathbf{U} - \mathbf{\bar{U}})\mathbf{\Sigma}_d  \|_F^2 = 2\sum_{i=1}^d \sigma_i (1 - \mathbf{u}_i^T \mathbf{\bar{u}}_i) = 2\sum_{i=1}^d \sigma_i (1 - \cos (\theta_i))$,
	where $\theta_i$ is the angle between column vectors $\mathbf{u}_i$ and $\mathbf{\bar{u}}_i$ (the columns of $\mathbf{U}$ and $\mathbf{\bar{U}}$, respectively), and $\sigma_i > 0$ are the diagonal elements of $\mathbf{\Sigma}_d$.}\hfill$\square$

Therefore, the proposed optimization objective function minimizes the weighted sum of the cosines of the angles between the original columns and their approximations. As such, low-rank or scaling issues cannot arise. The only potential problem is that different columns might have very different approximation errors (i.e., there could exist $i$ such that $\cos(\theta_i) \approx 1$ while there could be some $j$ for which $\cos(\theta_j) \ll 1$). This issue can be mitigated by increasing $g$ in \eqref{eq:approxu} or choosing carefully the indices $(i_k, j_k)$ where the proposed transformations operate.

Assume that we consider the operator norm as the approximation error, i.e., we want to minimize $\|\mathbf{U} - \mathbf{\bar{U}}\|_2$ in \eqref{eq:theoptimizationproblem_only_orthonormal}. We have the following two results.

\noindent \textbf{Remark 8. [The spectrum of the error matrix]} \textit{Consider $p = d$ and $\mathbf{\bar{\Sigma}}_d = \mathbf{\Sigma}_d = \mathbf{I}_d$ in \eqref{eq:theoptimizationproblem_only_orthonormal} then for any orthonormal $\mathbf{U}$ and $\mathbf{\bar{U}}$ the error matrix $\mathbf{U} - \mathbf{\bar{U}}$ is normal and has all its eigenvalues on a circle of radius one centered at $(1,0)$ in the complex plane. As such, we have that $\| \mathbf{U} - \mathbf{\bar{U}} \|_2 \leq 2$.}\hfill$\square$

\noindent \textbf{Theorem 5. [A bound on the operator norm of the error matrix]} \textit{Consider $p = d$ and $\mathbf{\bar{\Sigma}}_d = \mathbf{\Sigma}_d = \mathbf{I}_d$ in \eqref{eq:theoptimizationproblem_only_orthonormal} and that $\mathbf{u}^T_i \mathbf{\bar{u}}_i \geq 0$ for all $i=1,\dots,d$, then the operator norm of the error matrix obeys $\| \mathbf{U} - \mathbf{\bar{U}} \|_2 \leq 1-\epsilon_\text{min} + \sqrt{(d-1) (1-\epsilon_\text{min}^2) }$ where $\epsilon_\text{min}  = \underset{i}{\min} \ \mathbf{u}^T_i \mathbf{\bar{u}}_i = \underset{i}{\min} \ \cos(\theta_i)$. The bound in Remark 8 is met when $\epsilon_{\min} \geq (d-2)/d$.}\hfill$\blacksquare$

Remark 8 describes the full spectrum of the error matrix. It is interesting to notice that the upper bound $\| \mathbf{U} - \mathbf{\bar{U}} \|_2 \leq 2$ coincides with the expectation result from \cite{CollinsMale2011}: as $d \to \infty$ if both $\mathbf{U}$ and $\mathbf{\bar{U}}$ are chosen uniformly at random with Haar measure then almost surely $\|\mathbf{U} + \mathbf{\bar{U}} \| \to 2$. Theorem 5 describes an upper bound on the operator norm of the error matrix. The bound is tight only for very high values of $\epsilon_{\min}$ and therefore is meant to give a qualitative measure of the approximation. As in Remark 7, the key quantity is $\cos(\theta_i)$ but now the bound depends on the worst approximation: while the Frobenius norm objective function minimizes the sum of the pairwise distances between the columns of $\mathbf{U}$ and $\mathbf{\bar{U}}$, when using the operator norm the approximation accuracy depends on the largest distance between the same pairwise columns. The assumption that $\mathbf{u}^T_i \mathbf{\bar{u}}_i \geq 0$ is not restrictive at all as we expect $\mathbf{U}^T \mathbf{\bar{U}}$ to be diagonally dominant with positive elements on the diagonal (see Remark 5 and the discussion after Remark 6). We note that in this context the Frobenius norm approach could be used to minimize a proxy for the operator norm. As in Remark 7 we have the choice of $\sigma_i$ we could update these values with each iteration of the proposed algorithm such that $\sigma_i^{\text{(new)}} \leftarrow \frac{\mathbf{u}_{i_{\max}}^T \mathbf{\bar{u}}_{i_{\max}} }{\mathbf{u}_i^T \mathbf{\bar{u}}_i}$ where $i_{\max} = \underset{i}{\arg \max}\ \mathbf{u}_i^T \mathbf{\bar{u}}_i$. This choice will encourage the algorithm to improve upon the worst pairwise column approximation (i.e., the highest $\mathbf{u}_i^T \mathbf{\bar{u}}_i$). This approach would be in the spirit of an iteratively reweighted least squares algorithm such as \cite{doi:10.1002/cpa.20303}.

Finally, the last measure of accuracy we consider is the angle between two subspaces as described by \cite{10.2307/2005662}. This value can be non-zero only when $p<d$, a case which is interesting for PCA projections and which we detail in the next section. Following \cite{doi:10.1137/S1064827500377332}, given $\mathbf{U}_p$ and $\mathbf{\bar{U}}_p$, whose columns span two $d$-dimensional subspaces of size $p$, we compute the angles between the two subspaces as
\begin{equation}
\beta_i = \arccos (\tau_i),\ i=1,\dots,p,
\end{equation}
where $0 \leq \tau_i \leq 1$ are the singular values of $\mathbf{U}_p^T \mathbf{\bar{U}}_p$ and where $0\leq \beta_1 \leq \dots \leq \beta_p \leq \pi/2$. We take the principal angle to be the largest angle above, i.e., $\beta = \beta_p = \arccos(\tau_p)$.

%\noindent \textbf{Remark 9. [Connection between the operator norm and subspace distance errors]} \textit{Assuming that $\mathbf{U}_p^T \mathbf{\bar{U}}_p \approx \mathbf{I}_p$ we have that $\beta_i \approx 1-\tau_i$.}\hfill$\square$

Finally, we would like to note that all the approximation errors we have previously discussed measure (in different ways) how well $\mathbf{U}_p^T \mathbf{\bar{U}}_p$ approaches the identity matrix.

\section{Application: fast PCA projections}

Consider a  training set of $d$-dimensional points $\{\mathbf{x}_i \}_{i=1}^N$ and  the $d \times N$ matrix $\mathbf{X} = \begin{bmatrix} \mathbf{x}_1 & \dots & \mathbf{x}_N \end{bmatrix}$. Given $1 \leq p < d$, PCA provides the optimal $p$-dimensional projection that minimally distorts, on average, the data points. The projection is given by the  eigenvectors of the $p$ largest eigenvalues of $\mathbf{X} \mathbf{X}^T$ or, equivalently, the left singular vectors of the $p$ largest singular values of $\mathbf{X}$, that is
\begin{equation}
\mathbf{XX}^T \approx \mathbf{U}_p (\mathbf{\Sigma}_p \mathbf{\Sigma}_p^T ) \mathbf{U}_p^T \text{ and } \mathbf{X} \approx \mathbf{U}_p \mathbf{\Sigma}_p \mathbf{V}_p^T.
\label{eq:eigandsvd}
\end{equation}
Given the above decompositions we can approximate $\mathbf{X} $ by $\mathbf{\bar{X}} = \mathbf{\bar{U}} \mathbf{\bar{\Sigma}}_p \mathbf{V}_p^T$, i.e., we keep $\mathbf{V}_p$ but we modify the principal components and their singular values, such that we minimize the error given by
\begin{equation}
\| \mathbf{U}_p \mathbf{\Sigma}_p \mathbf{V}_p^T -  \mathbf{\bar{U}} \mathbf{\bar{\Sigma}}_p \mathbf{V}_p^T \|_F^2 = \| \mathbf{U}_p \mathbf{\Sigma}_p -  \mathbf{\bar{U}} \mathbf{\bar{\Sigma}}_p \|_F^2,
\label{eq:Visnotneeded}
\end{equation}
where $\mathbf{U}_p$ is $d \times p$, the diagonal matrix $\mathbf{\Sigma}_p$ is $p \times p$, $\mathbf{\bar{U}}$ is of size $d \times d$, $\mathbf{\bar{\Sigma}}_p$ is $d \times p$ and is zero except for its main $p \times p$ diagonal. In this paper, we work with $\mathbf{X}$, as opposed to $\mathbf{XX}^T$, to keep the relationship with $\mathbf{U}_p$ linear, rather than quadratic. In the context of applying our approach to PCA, we use our decomposition on the principal components $\mathbf{U}_p$ which we assume are already calculated together with the associated singular values $\mathbf{\Sigma}_p$ which we may use as weights in \eqref{eq:theoptimizationproblem_only_orthonormal}.\\
   %we propose a two step procedure \textit{Step 1}: given the dataset $\mathbf{X}$, get the optimal PCA components $\mathbf{U}_p$ and $\mathbf{\Sigma}_p$ by your favorite algorithm; \textit{Step 2}: given $\mathbf{U}_p$, $\mathbf{\Sigma}_p$ and a computational budget $g \ll d^2$, similarly to \eqref{eq:theoptimizationproblem_only_orthonormal}, we approximately solve \eqref{eq:theoptimizationproblem_only_orthonormal} for $\mathbf{U}_p \mathbf{\Sigma}_p$.
Note that in \eqref{eq:Visnotneeded}, $\mathbf{V}_p$, which has size $N$, is not necessary and that the two-step procedure is equivalent to computing the projections $\mathbf{\bar{U}}$ directly from $\mathbf{X}$. Also note that based on \eqref{eq:eigandsvd}, we could factor $\mathbf{X} \approx \mathbf{\bar{U}} \mathbf{\Sigma} \mathbf{\bar{V}}^T$ where $\mathbf{\bar{U}}$ and $\mathbf{\bar{V}}^T$ are approximations in $\mathcal{F}_{O(d \log d)}$ and $\mathcal{F}_{O(N \log N)}$, respectively. The difficulty here is the dependency of $\mathbf{\bar{V}}$ on $N \gg d$ which would require a large running time.\\% can be very large and therefore finding $\mathbf{\bar{V}}$ would require a large running time. It is for this reason that we choose to find $\mathbf{U}$ or $\mathbf{U}_p$ first by some well-established algorithms and then perform our approximation $\mathbf{\bar{U}}$ or $\mathbf{\bar{U}}_p$ depending only on dimensions $d$ and $p$.
With $\mathbf{\bar{U}}$ fixed, for $\mathbf{\bar{\Sigma}}_p$ we have several strategies: i) set it to the identity, i.e., flatten the spectrum; ii) keep it to the original singular values $\mathbf{\Sigma}_p$; or iii) continuously update it to minimize the Frobenius norm, i.e., get the new ``singular values'' that are optimal with the approximation $\mathbf{\bar{U}}$. The first approach favors the accurate reconstruction of all components while the other approach favors mostly the few leading components only (depending on the decay rate of the corresponding singular values).\\%The good choice depends on the application at hand.
\subsection{Comparison with the symmetric diagonalization by Givens rotations approach} Because of the locally optimal way the Givens transformations are chosen (see Theorem 1), our proposed factorization algorithm is computationally slower than the Jacobi diagonalization process which chooses the Givens rotations on indices $(i,j)$ corresponding to the largest off-diagonal entry of the covariance matrix. Furthermore, the Jacobi decomposition uses each rotation to zero the largest absolute value off-diagonal entry and because of this sub-optimal choice needs $O(d^2 \log d)$ Givens rotations \cite{ManySweeps} to complete de diagonalization (more than the $\frac{d(d-1)}{2}$ needed to fully represent the orthonormal group).\\
 In all other aspects, our approach provides advantages over the Jacobi approach: i) we define a clear objective function that we locally optimize exactly; ii) it is known that the Jacobi process converges slowly when the number of rotations is low \cite{JacobiIsLInearEarlyOn}, which is exactly the practically relevant scenario we have, i.e., $g \ll d^2$; iii) with the same computational complexity, i.e., $g$ terms in the factorization, our proposed approach is always more accurate since we include as a special case the Givens rotations.
	\subsection{Comparison with  structured matrix factorization} Our approach requires the explicit availability of the orthonormal principal directions $\mathbf{U}_p$. Previous methods that factor using only Givens rotations are not applied directly on an orthogonal matrix. These methods rely on receiving as input a symmetric object (e.g., $\mathbf{XX}^T$) and then using a variant of Jacobi iterations for matrix diagonalization \cite{lemagoarou:hal-01416110} or multiresolution factorizations \cite{Kondor2014MMF} to find the good rotations that approximate the orthonormal eigenspace. Applying Givens transformations directly to $\mathbf{X}$ on the left, i.e., $\mathbf{G}_{ij} \mathbf{X}$, cannot lead to the computation of the PCA projections $\mathbf{U}$ but only to the polar decomposition. On the other hand, when applying Givens rotations on both sides of the covariance matrix, i.e., $\mathbf{G}_{ij} \mathbf{XX}^T \mathbf{G}_{ij}^T$, then the right eigenspace $\mathbf{V}$ cancels out in the product \eqref{eq:eigandsvd} and we are able to directly recover $\mathbf{U}$ (but we need $\mathbf{XX}^T$ explicitly). Finally, note that the diagonalization process approximates the full eigenspace $\mathbf{U}$ and cannot separate from the start the $p$ principal components $\mathbf{U}_p$ because they are solving the following problem
	\begin{equation}
	\underset{\mathbf{\bar{U}},\ \mathbf{\bar{\Lambda}}}{\text{minimize}}\ \|  \mathbf{XX}^T - \mathbf{\bar{U}}\mathbf{\bar{\Lambda}} \mathbf{\bar{U}} \|_F^2 \text{ subject to } \mathbf{\bar{U}} \in \mathcal{F}_g.
	\end{equation}
	This formulation is useful to approximate the whole symmetric  matrix $\mathbf{XX}^T$ (or $\mathbf{U}$), but not necessarily the $p$ principal eigenspace $\mathbf{U}_p$. To get these we would need to complete the diagonalization process, find the $p$ largest entries on the diagonal of $\mathbf{\bar{\Lambda}}$ and then work backward  to identify  the rotations that contributed  diagonalizing those largest elements. This  procedure would be prohibitively  expensive. Previous work, e.g. \cite{Kondor2014MMF, Treelets}, deals with approximating $\mathbf{XX}^T$ rather than computing PCA. In the same line of work, the one-sided Jacobi algorithm for SVD \cite{doi:10.1137/05063920X} can also be applied. 

\section{Experimental results}

Given the rather pessimistic guarantees, we tackle problems: how well does Algorithm 1 recover random orthogonal matrices and principal components such that we benefit from the computational gains but do not significantly impact the approximation/classification accuracy. Source code available.%\footnote{Matlab source code: https://gofile.io/?c=wqf3L9}
\footnote{https://github.com/cristian-rusu-research/fast-orthonormal-approximation}

\subsection{Synthetic experiments}

\begin{figure}[!tbp]
	\centering
	\begin{minipage}[t]{0.495\textwidth}
		\includegraphics[trim = 18 5 30 22, clip, width=0.49\textwidth]{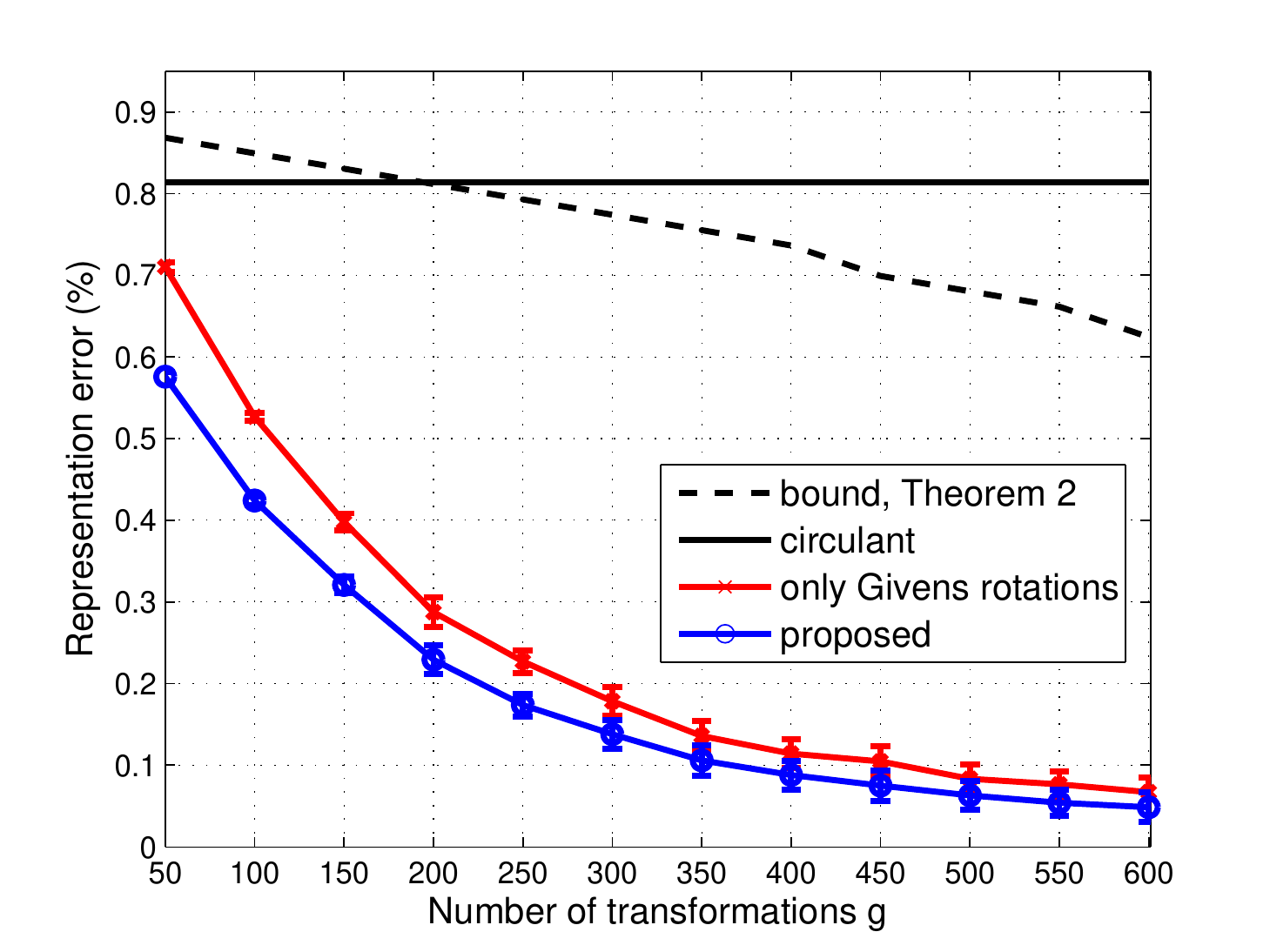}
		\includegraphics[trim = 18 5 28 22, clip, width=0.49\textwidth]{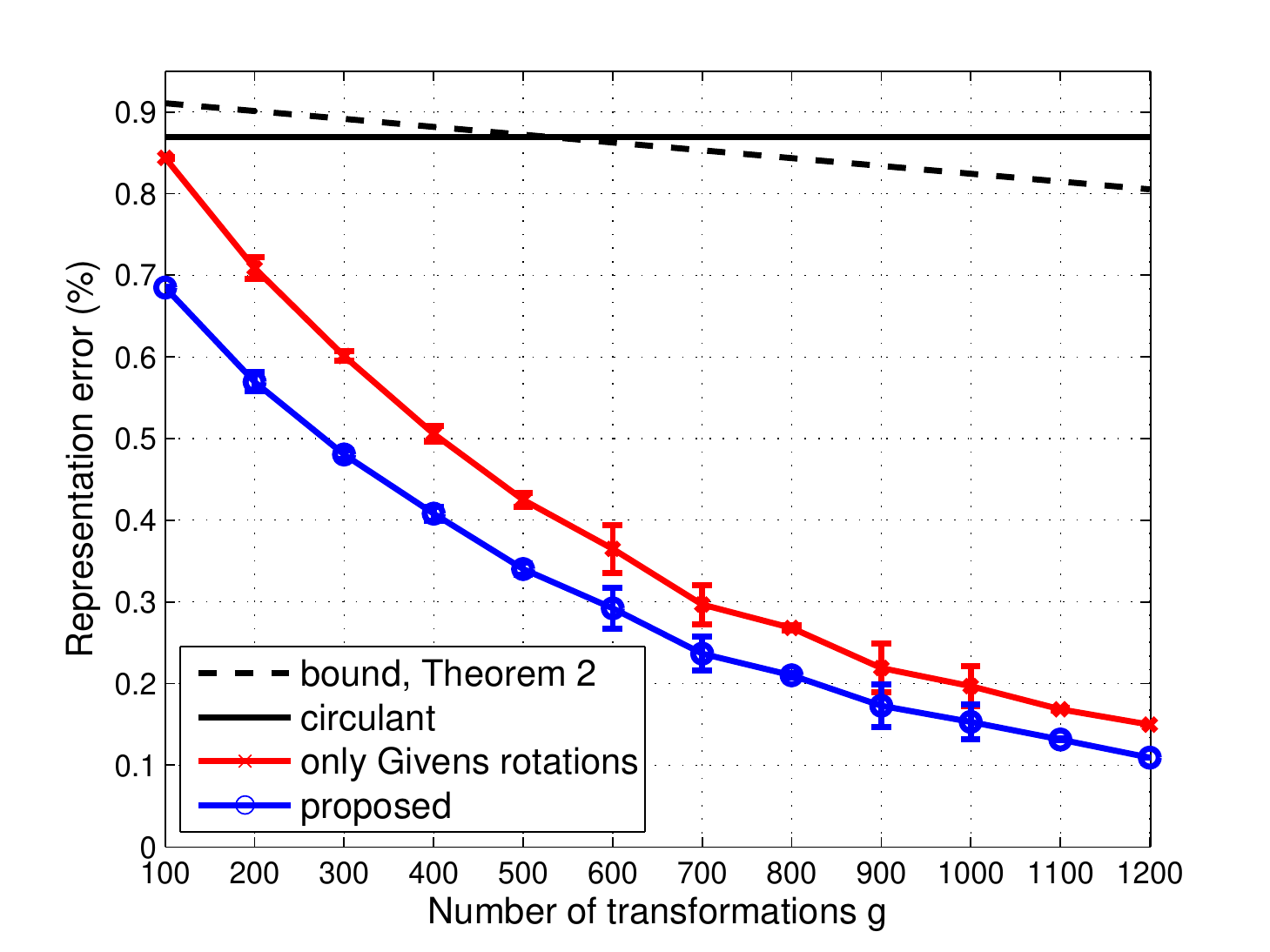}
		\caption{Average approximation errors and standard deviations over 100 realizations of random orthonormal matrices of size $d = 50$ (left) and $d = 100$ (right). For reference we show the bound developed in Theorem 3, the approximation accuracy of the circulant \cite{CirculantDRoperator} (Toeplitz performed just marginally better) and that of the using the factorization \eqref{eq:approxu} but allowing only Givens rotations. As expected, performance degrades with large $d$.}
		\label{fig:random}
	\end{minipage}
	\hfill
	\begin{minipage}[t]{0.495\textwidth}
		\includegraphics[trim = 18 5 28 20, clip, width=0.49\textwidth]{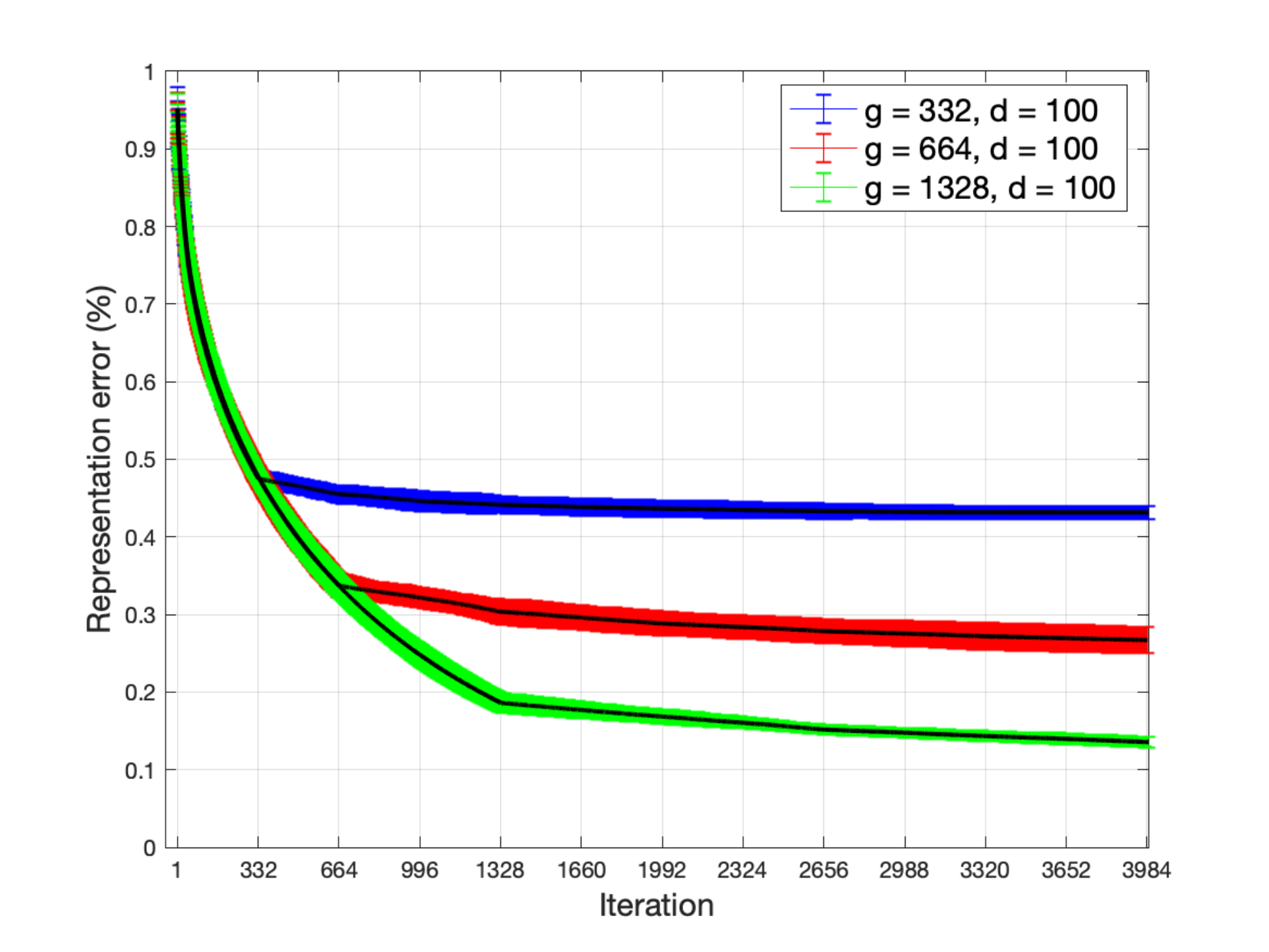}
		\includegraphics[trim = 18 5 30 20, clip, width=0.49\textwidth]{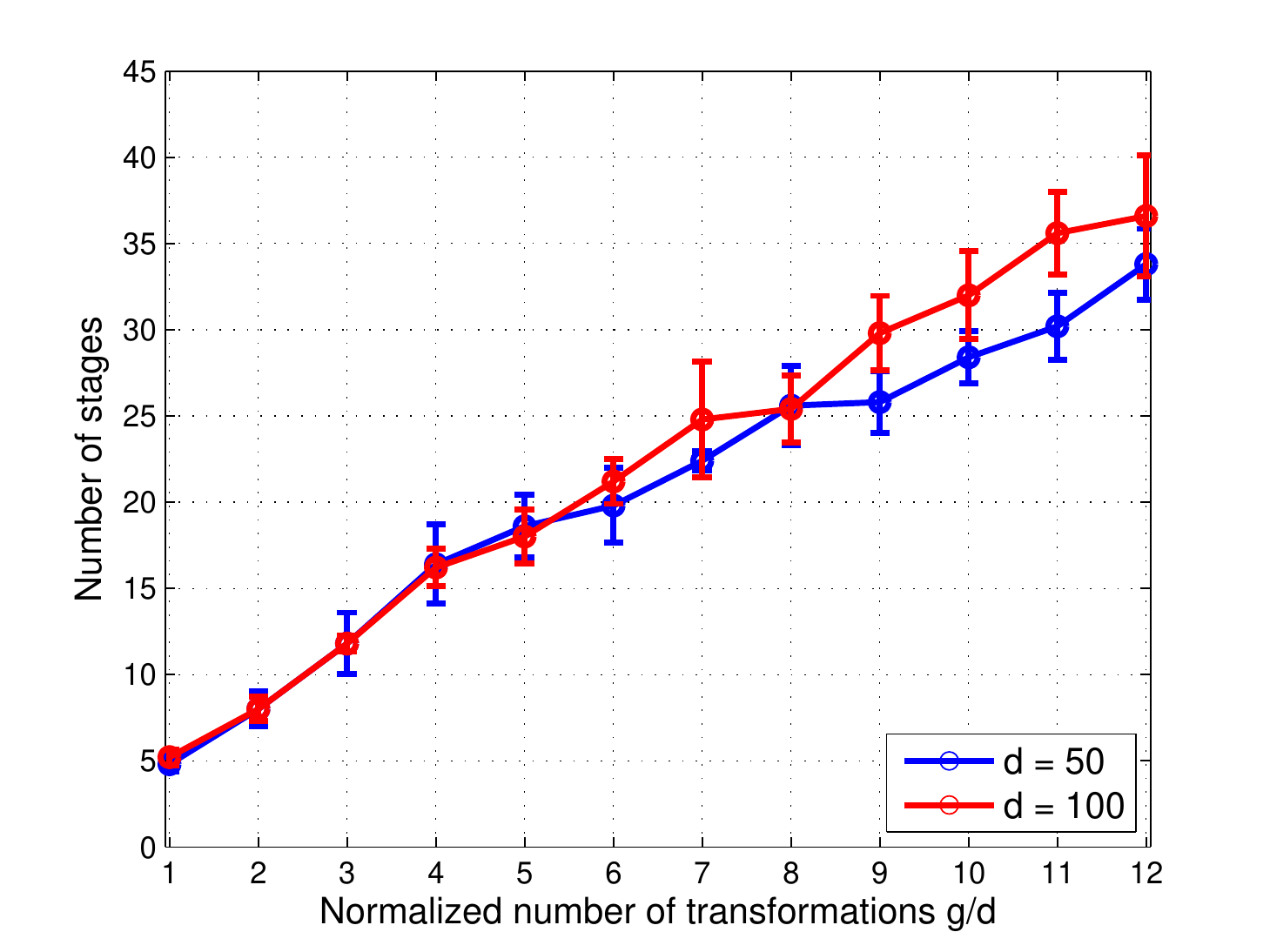}
		\caption{Left: for $d=100$ and $g \in \{   332, 664, 1328 \}$ we show the evolution (mean and standard deviation) of the objective function with the number of iterations. In each case, the first $g$ iterations are the initialization process. Right: for $d \in \{50, 100\}$ we show the number of stages in each transformation we  learrn with the proposed  algorithm as a function of $g$. For both plots results are averaged over 100 realizations.}
		\label{fig:random_evo_stages}
	\end{minipage}
\end{figure}

For fixed $d$ we generate random orthonormal matrices from the Haar measure \cite{HowToRandomUnitary}. %as follows: i) we generate a random $d \times d$ matrix with entries drawn i.i.d. standard Gaussian; ii) we apply the QR algorithm and we keep the orthogonal factor $\mathbf{Q}$; finally iii) we build $\mathbf{U} = \mathbf{DQ}$ where $\mathbf{D}$ is a diagonal matrix with entries $\pm 1$ such that $\mathbf{U}$ has Haar measure  \cite{HowToRandomUnitary}.
Figure \ref{fig:random} shows the representation error $(2d)^{-1} \| \mathbf{U} - \mathbf{\bar{U}} \|_F^2$ for the proposed method. The plot shows that allowing for the Givens transformations $\mathbf{\tilde{G}}_{i j}$ in \eqref{eq:theG} brings a 17\% relative benefit as compared to using only the Givens rotations while, for the same $g$, the computational complexity is the same. The circulant approximation performs worst because it has the lowest number of degrees of freedom, only $d$ (computationally, it is comparable with the Givens and proposed approaches for $g = 100$). Lastly, we can observe that the bound is very pessimistic, especially for these values of $g$. In Figure \ref{fig:random_evo_stages} (left) we show for fixed number of transformations $g$ the progress that the proposed algorithm makes with each iteration. It is interesting to observe that the initialization steps (first $g$ steps) significantly decrease the approximation error while the other step make only moderate improvements. This indicates that convergence is slow and might take a large number of iterations with little progress made by the latter steps. Also in Figure \ref{fig:random_evo_stages} (right) we show the number of stages in each transformation. A stage is a set of extended orthogonal Givens transformations that can be applied in parallel, i.e., they do not share any indices $(i_k, j_k)$ among them. This is important from an implementation perspective as parallel processing can be exploited to speedup the transformations.

\subsection{MNIST digits and fashion}
\begin{figure}[!tbp]
	\centering
	\begin{minipage}[t]{0.495\textwidth}
		\includegraphics[trim = 18 5 26 20, clip, width=0.49\textwidth]{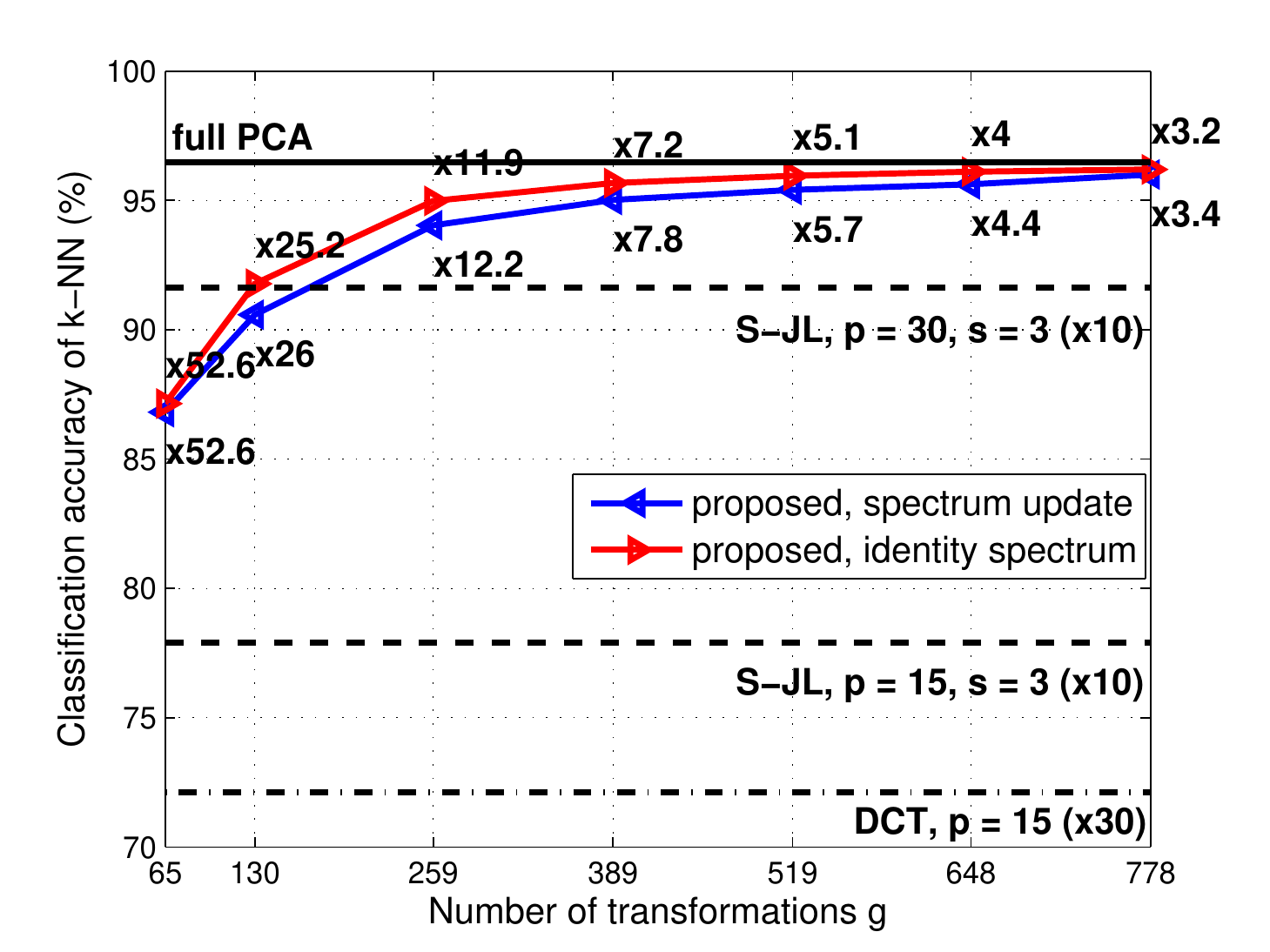}
		\includegraphics[trim = 18 5 23 20, clip, width=0.49\textwidth]{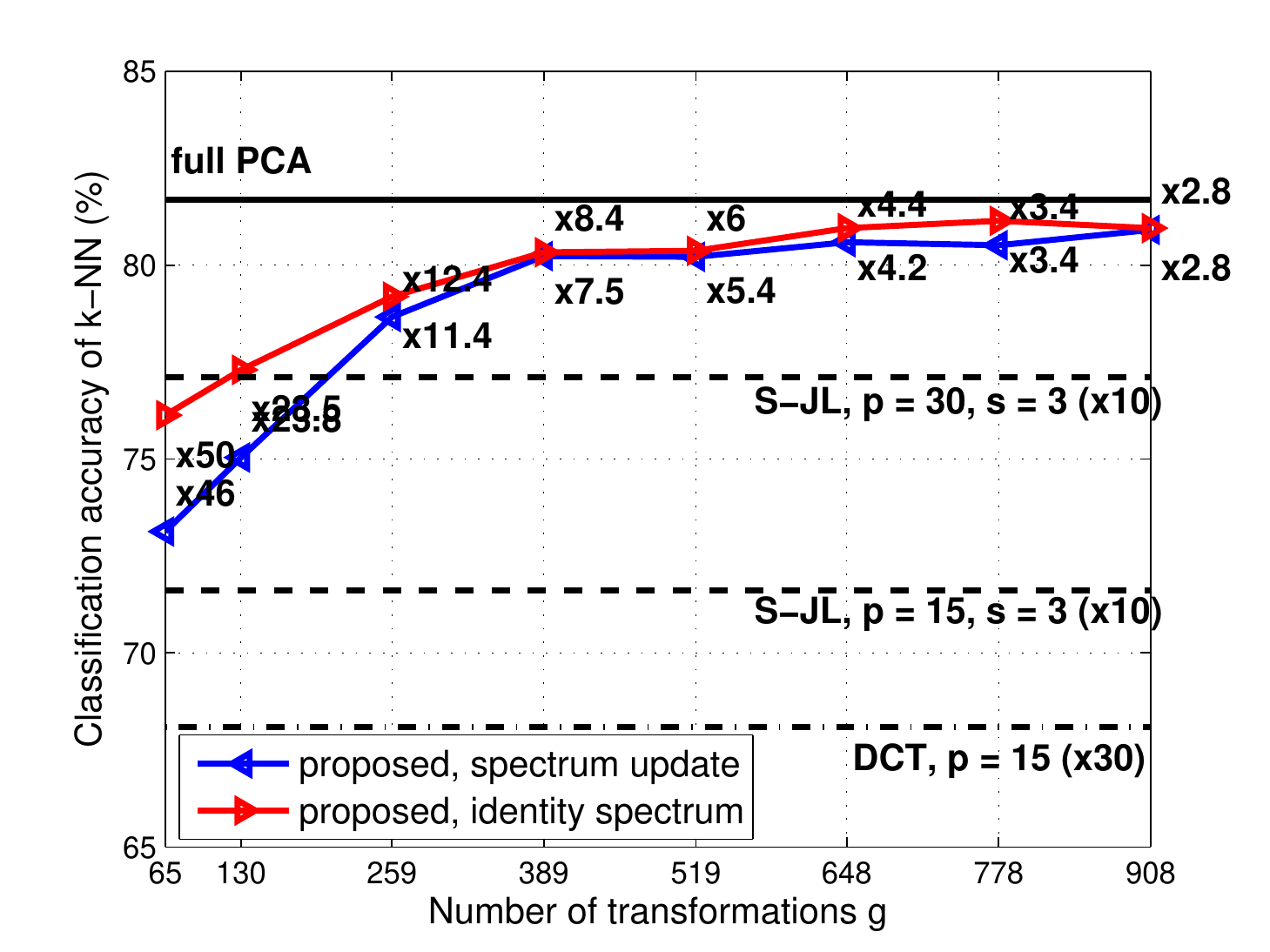}
		\caption{Classification accuracy obtained by the k-NN algorithm for the MNIST digits (left) and fashion (right) datasets as a function of the complexity of the proposed projections. Dimensionality reduction was done with $p = 15$ principal components. The bold text represents the speedup (FLOPS) compared to the cost of projecting with the unstructured, optimal, PCA components which takes $2pd$ operations. For the sparser JL the variable $s$ is the number of non-zeros in each column.}
		\label{fig:mnistdetails}
	\end{minipage}
	\hfill
	\begin{minipage}[t]{0.495\textwidth}
		\includegraphics[trim = 18 5 28 20, clip, width=0.49\textwidth]{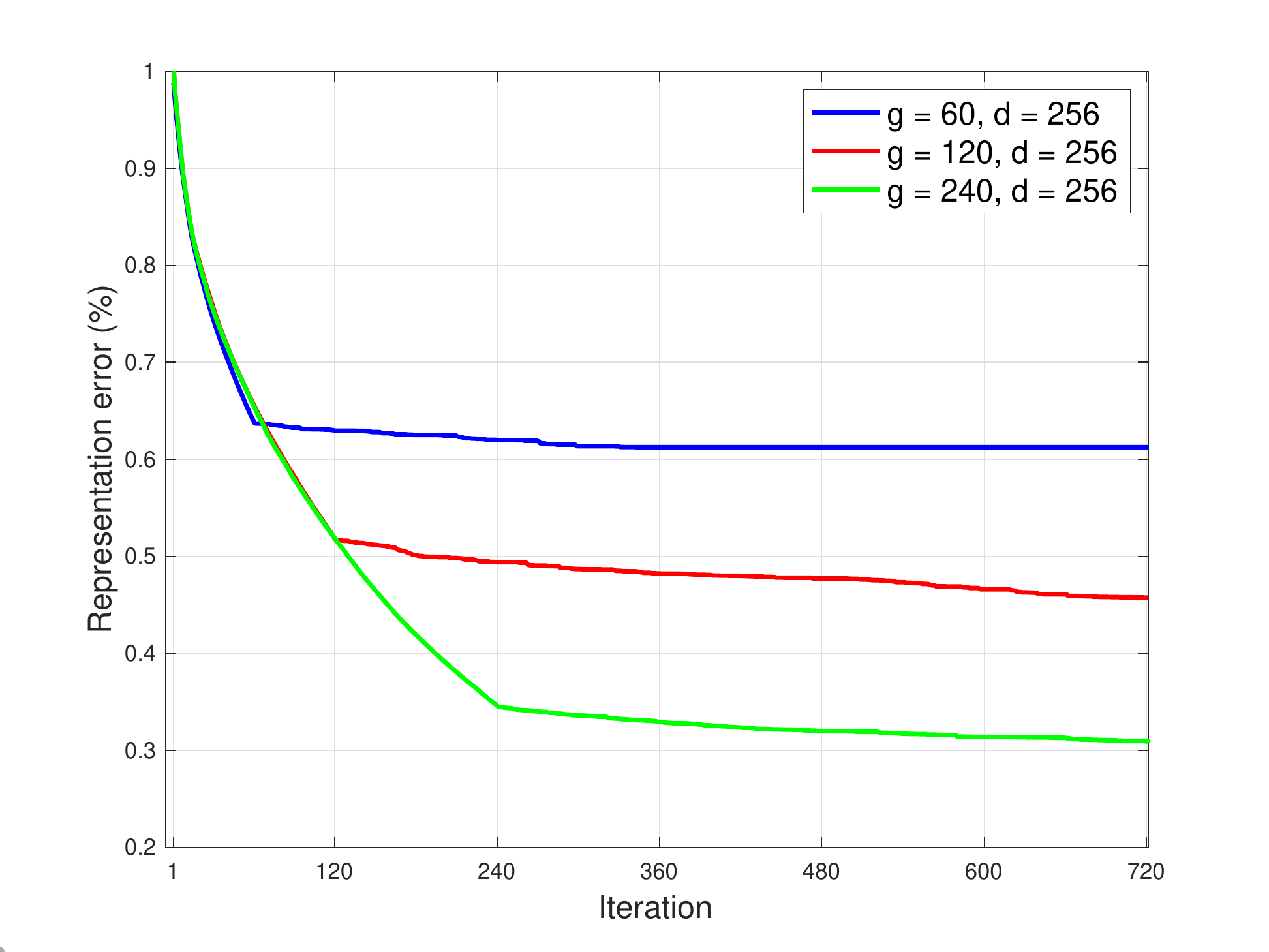}
		\includegraphics[trim = 18 5 30 20, clip, width=0.49\textwidth]{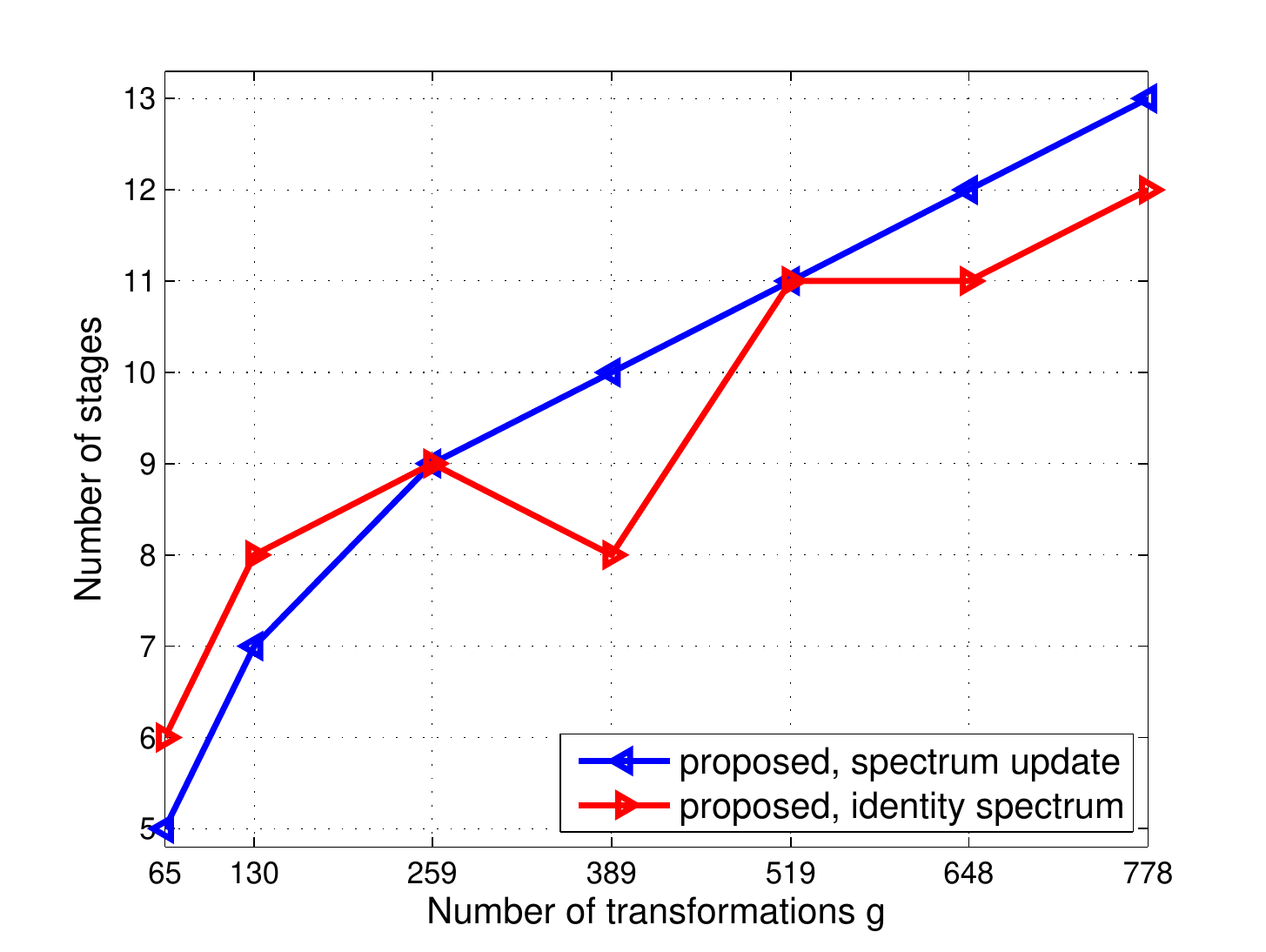}
		\caption{For the USPS dataset ($d  =  256$) we have, left:  and $g \in \{  60, 120, 240 \}$ we show the evolution (mean and standard deviation) of the objective function with the number of iterations -- again, in  each case, the first $g$ iterrations are the initialization steps; right: the number of stages in each transformation created by the proposed method as a function of the number of extended orthogonal Givens transformations $g$.}
		\label{fig:usps_evo_stages}
	\end{minipage}
\end{figure}

We now turn to a classification problem. We use the MNIST digits
%\footnote{http://yann.lecun.com/exdb/mnist}
 and fashion
 %\footnote{https://github.com/zalandoresearch/fashion-mnist}
  datasets. The points have size $d = 400$ (we trimmed the bordering whitespace) and we have $N = 6\times 10^4$ training and $N_\text{test} = 10^4$ test points. In all cases, we use the k-nearest neighbors (k-NN) algorithm with $k = 10$, and we are looking to correctly classify the test points. Before k-NN we apply PCA and our proposed method. Results are shown in Figure \ref{fig:mnistdetails}. We deploy two variants of the proposed method: approximate the principal components as if they had equal importance and approximate the principal components while simultaneously also updating estimates of the singular values.\\
For comparison, we also show the sparse JL \cite{SparserJL}. In this case, the target dimension is $p \in \{15, 30\}$ while the random transformation of size $p \times d$ only has three non-zero entries per column. More non-zeros did not have any significant effect on the classification accuracy while increasing (doubling, in this case) the target dimension $p$ increases the classification accuracy by 10\%. The results reported in the plots are averages for 100 realizations and the standard deviation is below 1\%. Of course, the significant advantage of the JL approach over PCA is that no training is needed. The disadvantage is that if we choose greedily the target projection dimension, i.e., low $p$, the accuracy degrades significantly for JL. Results are identical for PCA when $p$ is $15$ of $30$.\\
Since we are dealing with image data, we also project using the discrete cosine transform (DCT). For the digits dataset, the performance is poor but, surprisingly, for the fashion dataset, this approach is competitive given the large speedup it reports (we used the sparse fast Fourier transform \cite{Indyk:2014:SSF:2634074.2634110} as we want only the largest $p$ components). We have also performed the projection by fast wavelet transforms with `haar' and several of the Daubechies `dbx' filters but the results were always similar to that of the DCT. For clarity of exposition, we did not add these results to the figures.\\
Our proposed methods report a clear trade-off between the classification accuracy and the numerical complexity of the projections. If we insist on an accuracy level close (within 1--2\%) of the full PCA then the speedup is only about x3. Reasonable accuracy is obtained for a speedup of x4--x5 after which the results degrade quickly. For $p=15$ better performance seems impossible via randomization. In this figure, the speedup is measured in terms of the number of operations (FLOPS).\\
Finally, in Figure \ref{fig:gt_vs_spca} (left) we compare our proposed methods against the sparse PCA on MNIST digits. sPCA performs exceptionally well in terms of the classification accuracy given the computational budget (a similar result is replicated for MNIST fashion). On other datasets where the principal components capture some global features (not local like in our example) we expect this performance to degrade. The training time of sPCA exceeds by 60\% the running time of PCA plus that of our method. We used the implementation of \citet{MINIMAXOPTIMALSPCA}.
\begin{table}[t]
	\caption{Average classification accuracies for k-NN when using PCA projections and our approximations without spectrum update, for various datasets. We also show the speedup (FLOPS and actual running time) and the number of features selected in the calculations as a proportion out of the total $d$ (see also Remark 2).  Results are averaged over 100 random realizations (train/test splits).}
	\label{tb:sample-table}
	\vskip 0.15in
	\begin{center}
		\begin{small}
			\begin{sc}
				\begin{tabular}{lcccccr}
					\toprule
					DATASET & FULL PCA & \multicolumn{4}{c}{PROPOSED ALGORITHM}  \\
					& accuracy & accuracy  & speedup(FLOPS) & speedup(TIME)  & selection \\
					\midrule
					PENDIGITS    & 95 $\pm$ 0.7& 91 $\pm$ 2.1 & $\times$1.6 & $\times$1.1 & 1 \\ % file is: penset p = 7 m = 28
					ISOLET & 92 $\pm$ 0.4 & 90 $\pm$ 1.0 & $\times$12 & $\times$10.1 & 1 \\ % file is: isolet p = 150 m = 2781
					USPS    & 95 $\pm$ 1.2& 94 $\pm$ 0.8 & $\times$7.7 &  $\times$4.7 & 0.64 \\ % file is: usps p = 15 m = 480
					UCI    & 90 $\pm$ 1.9& 87 $\pm$ 1.5 & $\times$2.5 & $\times$1.6  & 0.72 \\ % file is: uci p = 6 m = 72
					20NEWS     & 80 $\pm$ 3.1 & 77 $\pm$ 2.1 & $\times$3.1  & $\times$2.5 & 0.3 \\ % file is: 20news p = 200 m = 3480
					EMNIST digits & 97 $\pm$ 2.5& 95 $\pm$ 1.8 & $\times$13 & $\times$11.3 & 0.37 \\ % file is: emnist digits p = 15 m = 288
					MNIST 8m     & 96 $\pm$ 2.0& 94 $\pm$ 0.9 & $\times$15 & $\times$13.7 & 0.28  \\ % files is: mnist8m p = 15 m = 865
					\bottomrule
				\end{tabular}
			\end{sc}
		\end{small}
	\end{center}
	\vskip -0.1in
\end{table}
Because of ideas in Remark 1, we perform only the finally useful calculations and further reduce the computational cost on average by one third (these are accounted for already in the numbers in the plots).

\subsection{Experiments on other datasets}

\begin{figure}[!tbp]
	\centering
	\begin{minipage}[t]{0.495\textwidth}
		\includegraphics[trim = 18 5 28 20, clip, width=0.49\textwidth]{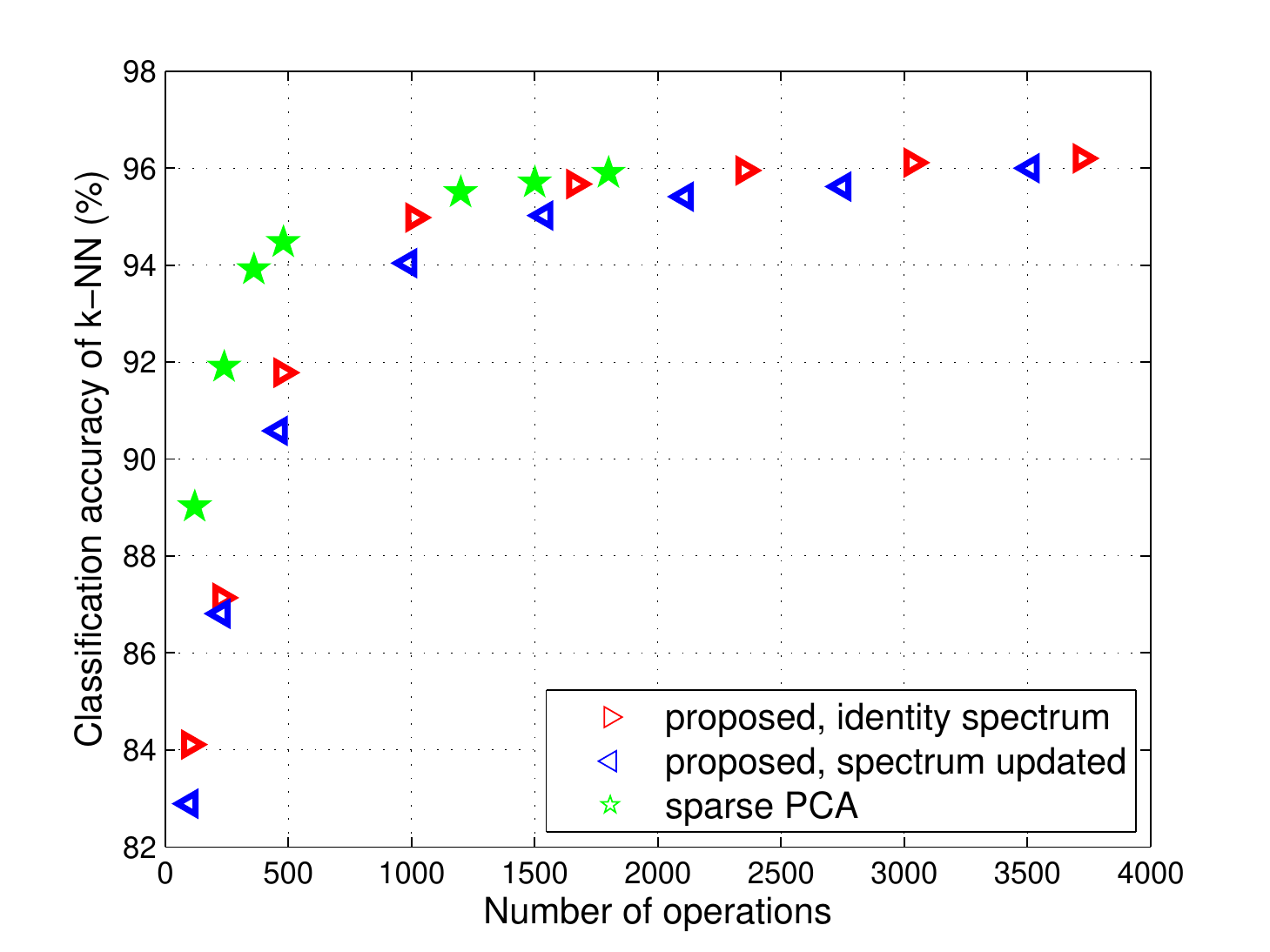}
		\includegraphics[trim = 18 5 30 20, clip, width=0.49\textwidth]{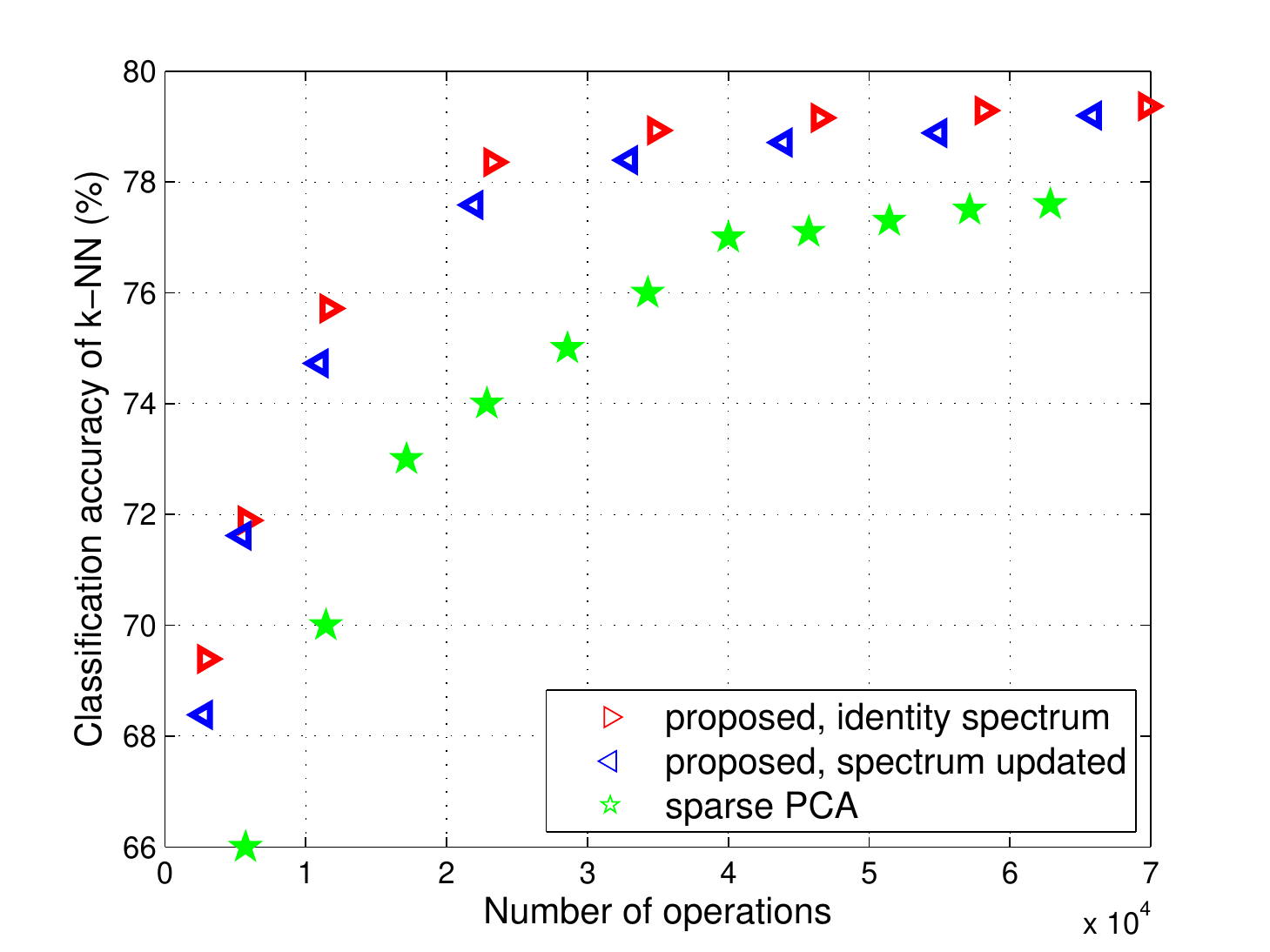}
		\caption{Classification accuracy versus number of operations on the MNIST digits (left) and 20NEWS (right) datasets for our proposed methods and the sparse PCA method \cite{MINIMAXOPTIMALSPCA, SparsePCA}.}
		\label{fig:gt_vs_spca}
	\end{minipage}
	\hfill
	\begin{minipage}[t]{0.495\textwidth}
		\centering
		\includegraphics[trim = 18 15 30 20, clip, width=0.75\textwidth]{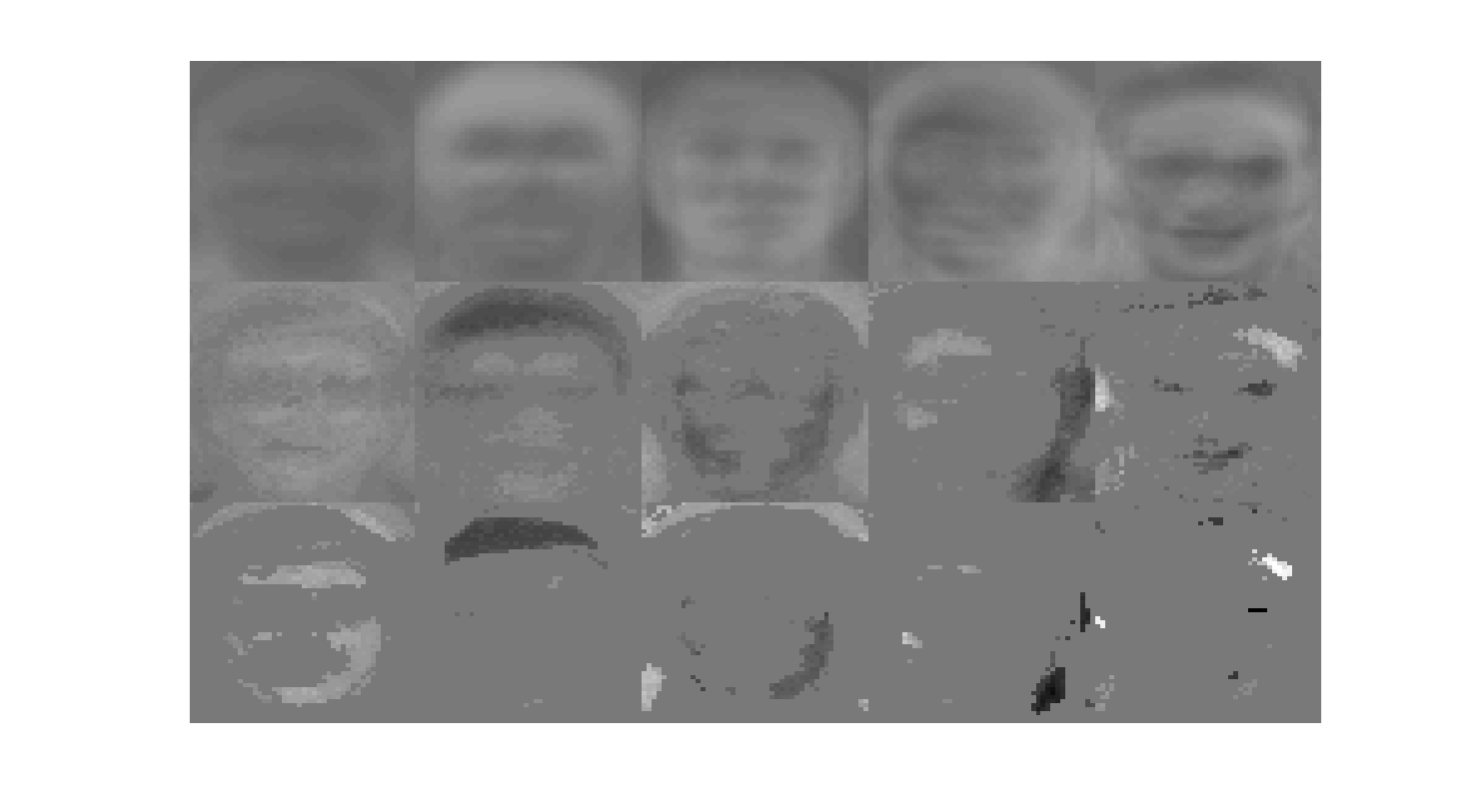}
		\caption{A few eigenfaces obtained by the optimal PCA (top) and by our proposed method with $g = 3059$ and $g = 1020$ (middle and bottom, respectively). The projection speedup (FLOPS) is x4.1 and x13.9, respectively.}
		\label{fig:eigenfaces}
	\end{minipage}
\end{figure}

The 20-newsgroups
%\footnote{http://qwone.com/$\sim$jason/20Newsgroups}
 dataset consists of $18827$ articles from 20 newsgroups (approximately 1000 per class). The data set was tokenized using the rainbow package (www.cs.cmu.edu/~mccallum/bow/rainbow).
 %\footnote{www.cs.cmu.edu/~mccallum/bow/rainbow}
 Each article is represented by a word-count vector for the $d = 2\times10^4$ common words in the vocabulary. For this dataset we have $N_\text{test} = 5648$, $p = 200$, and as shown in Figure \ref{fig:gt_vs_spca} (right) in this case we outperform sparse PCA.\\
We also apply our algorithm to several other popular dataset from the literature: PENDIGITS (www.ics.uci.edu/$\sim$mlearn/MLRepository.html)
%\footnote{www.ics.uci.edu/$\sim$mlearn/MLRepository.html}
 with 10 classes and $d = 16$, $N = 7494$, $N_\text{test} = 3498$, $p = 4$; ISOLET (archive.ics.uci.edu/ml/datasets/isolet)
 %\footnote{archive.ics.uci.edu/ml/datasets/isolet}
  with 26 classes and $d = 617$, $N = 6238$, $N_\text{test} = 1559$, $p = 150$; USPS (github.com/darshanbagul/USPS\_Digit\_Classification)
  %\footnote{github.com/darshanbagul/USPS\_Digit\_Classification\newline}
   with 10 classes and $d = 256$, $N = 7291$, $N_\text{test} = 2007$, $p = 12$; UCI (ftp.ics.uci.edu/pub/machine-learning-databases/optdigits)
   %\footnote{ftp.ics.uci.edu/pub/machine-learning-databases/optdigits}
    with 10 classes and $d = 64$, $N = 3823$, $N_\text{test} = 1797$, $p = 6$; EMNIST digits (nist.gov/itl/iad/image-group/emnist-dataset)
    %\footnote{nist.gov/itl/iad/image-group/emnist-dataset}
     with 10 classes and $d = 784$, $N = 24\times10^4$, $N_\text{test} = 4\times10^4$, $p = 15$; MNIST 8m (leon.bottou.org/papers/loosli-canu-bottou-2006)
     %\footnote{leon.bottou.org/papers/loosli-canu-bottou-2006}
      with 10 classes and $d = 784$, $N = 6.4 \times 10^6$, $N_\text{test} = 1.7 \times 10^6$, $p = 15$. All the results are shown in Table \ref{tb:sample-table}. Here we provide two measures for the speedup: the FLOPS (number of arithmetic operations, additions and multiplications) and the actual running time (in seconds). For the time speedup, we first computed $\mathbf{\bar{U}}$ and then implemented the matrix-vector multiplication and compared it to the generic matrix-vector multiplication with $\mathbf{U}$, both compiled in C using the gcc compiler and --O3 flag. Once computed by the proposed algorithm, the transformation could also be hardcoded and a further speedup improvement could be achieved. The running time speedups are slightly below the FLOPS speedups due to overhead in modern CPUs (fetching data, register loading etc.). Aside the computational benefits of the proposed transformations we note that while complex instruction set computing machines are closing the speedup gap in terms of running time the price is higher energy consumption and more complex execution pipelines/circuitry.\\
\begin{figure}[!tbp]
	\centering
	\includegraphics[trim = 110 0 132 12, clip, width=\textwidth]{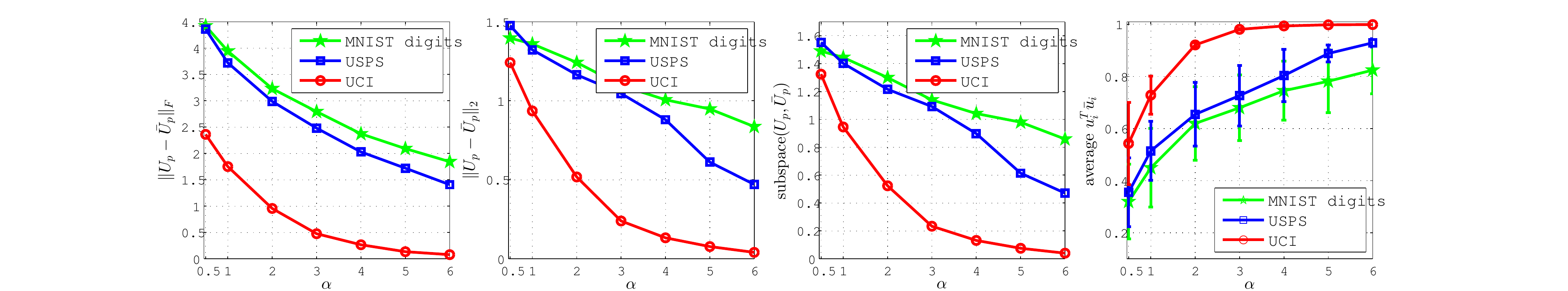}
	\caption{We show for three datasets used in Section 5.3, as a function of the number of Givens transformations $g$ in $\mathbf{\bar{U}}_p$ that goes like $g = \alpha p \log_2 d$, the Frobenius (left most) and operator norms errors (second from the left), angle between subspaces distance measured in radians (third from the left) and correlations (right most shows the average correlations but also minimum and maximum values). The results are averaged over 10 realizations (random training/testing data samples) but the variance is always below 0.05.}
	\label{fig:all_other_norms}
\end{figure}
In Figure \ref{fig:usps_evo_stages}, for the  USPS dataset, we provide insights into the behavior of the proposed algorithm with each iteration and the number of stages in each transformation that we construct. Results are similar to the ones shown for the random orthogonal approximation. We observe again that the initialization step is very efficient in reducing the approximation error (especially when compared to the iterative process that follows) and that the number of stages that have to be applied sequentially can further improve the running time when implementing the proposed transformations.\\
Finally, one of the most famous applications of PCA is in the field of computer vision for the problem of human face recognition. The eigenfaces \cite{Sirovich:87} approach was used successfully for face recognition and classification tasks. Here, we want to reproduce the famous eigenfaces by using the proposed methods. The original eigenfaces and their approximations (sparse eigenfaces \cite{SparseFaces}), with different $g$ and therefore different levels of detail, are shown in Figure \ref{fig:eigenfaces}. 

\subsection{Results on other approximation errors}

In Section 3.1 we have explored several other ways to measure the approximation accuracy of the proposed algorithm. We note that in the proposed method we keep the Frobenius norm objective function but we also measure the operator norm, the angle between subspaces distance and the average/minimum/maximum correlations. All results for three datasets are shown in Figure \ref{fig:all_other_norms}, the choices of parameters are the same as in Section 5.3. The first observation is that the approximation errors increase with larger $d$ (as already clear from Remark 4 and Theorem 3), i.e, best results are obtained for UCI ($d = 64$)  and the worst for MNIST ($d = 400$). Second, note that improving the Frobenius norm error we also improve all the other error measures. Thirdly, note the similar behavior between the operator norm and subspace distance errors. It is interesting to observe that the two plots seem to suggest experimentally that $\arccos(\sigma_{\min}( \mathbf{U}_p^T \mathbf{\bar{U}}_p )) \approx 1-\sigma_{\max}( \mathbf{I}_p - \mathbf{U}_p^T \mathbf{\bar{U}}_p )$ which hold over a large array of number of Givens transformations $g$. Lastly, the two right-most plots highlight the connection between the operator norm optimization and the maximization of the lowest coherence between columns and their approximations (as per Theorem 5).\\
Depending on the application at hand we can choose how to measure the approximation error.

\subsection{Details of the implementation}

\begin{table}[t]
	\caption{Speed-up achieved with the proposed approximations $\mathbf{\bar{U}}_p$ against a vanilla C implementation of dense matrix-vector multiplication and BLAS Level 2 functions, in both cases with no parallelism.}
	\label{}
	\vskip 0.15in
	\begin{center}
		\begin{small}
			\begin{sc}
				\begin{tabular}{lcccccr}
					\toprule
					DATASET: & ISOLET & USPS  & EMNIST digits & MNIST 8m  \\
					\midrule
					C language (TIME) & $\times$10.1 & $\times$4.7 & $\times$11.3 & $\times$13.7 & \\
					BLAS (TIME) & $\times$4.8 & $\times$2.5 & $\times$5.7 & $\times$6.8 & \\
					\bottomrule
				\end{tabular}
			\end{sc}
		\end{small}
	\end{center}
	\vskip -0.1in
\end{table}
The source code attached to this paper is written for the Matlab environment. Besides the implementation of Algorithm 1, we also provide the code to efficiently perform the matrix-vector multiplication with the proposed $\mathbf{\bar{U}}$ and $\mathbf{\bar{U}}_p$, respectively (also taking into account Remark 3). Still, due to the characteristics of Matlab, our implementation is not faster than the dense matrix-vector multiplication ``*''. To provide an appropriate comparison, we implement the matrix-vector multiplication with our structures in the compiled lower-level programming language C. We provide comparisons of this implementation against two scenarios: a vanilla matrix-vector multiplication (using the appropriate ordering of the loops to exploit locality) and the Basic Linear Algebra Subprograms (BLAS) Level 2 routines for matrix-vector multiplication (SGEMV). We use the single thread variant of BLAS as for the dimensions $d$ we consider the overhead of the parallel implementation is significant.\\
Applying the proposed transformations on batch features, i.e., matrix-matrix multiplications, would require careful application of the proposed transformations \eqref{eq:theG} to fully use the computing architecture and maximize performance.

\section{Conclusions}
This paper proposes a new matrix factorization algorithm for orthogonal matrices based on a class of structured matrices called extended orthogonal Givens transformations. We show that there is a trade-off between the computational complexity and accuracy of the approximations created by our approach. We apply our method to the approximation of a fixed number of principal components and show that, with a minor decrease in performance, we can reach significant computational benefits.\\
Future research directions include strengthening the theoretical guarantees since they  are way above what we observe experimentally. Further, it would be of interest to improve the  complexity of the proposed algorithm either by a parallel implementation or using randomization (e.g. computing a random subset of the $O(d^2)$ scores). As an immediate application, it would be interesting to apply our decomposition to the recently proposed unitary recurrent neural networks \cite{URNN}.

%\section*{Acknowledgements}
% We acknowledge the financial support of the AFOSR projects FA9550-17-1-0390 and BAA-AFRL-AFOSR-2016-0007 (European Office of Aerospace Research and Development), and the EU H2020-MSCA-RISE project NoMADS - DLV-777826. 

%\section*{Acknowledgements}

\section*{Acknowledgment}
This material is based upon work supported by the Center for Brains, Minds and Machines (CBMM), funded by NSF STC award CCF-1231216, and the Italian Institute of Technology. Part of this work has been carried out at the Machine Learning Genoa (MaLGa) center, Universita di Genova (IT). Lorenzo Rosasco acknowledges the financial support of the European Research Council (grant SLING 819789), the AFOSR projects FA9550-17-1-0390 and BAA-AFRL-AFOSR-2016-0007 (European Office of Aerospace Research and Development). Cristian Rusu acknowledges support by the Romanian Ministry of Education and Research, CNCS-UEFISCDI, project number PN-III-P1-1.1-TE-2019-1843, within PNCDI III.


%\bibliography{fast_pca}
%\bibliographystyle{icml2019}

%\section*{References}
{\small 
\bibliographystyle{plainnat}
\bibliography{arxiv_2019_2}

\begin{thebibliography}{52}
\providecommand{\natexlab}[1]{#1}
\providecommand{\url}[1]{\texttt{#1}}
\expandafter\ifx\csname urlstyle\endcsname\relax
  \providecommand{\doi}[1]{doi: #1}\else
  \providecommand{\doi}{doi: \begingroup \urlstyle{rm}\Url}\fi

\bibitem[Ailon and Chazelle(2006)]{FastJL}
N.~Ailon and B.~Chazelle.
\newblock Approximate nearest neighbors and the fast {J}ohnson-{L}indenstrauss
  transform.
\newblock In \emph{Proceedings of the 38th Annual ACM Symposium on Theory of
  Computing}, pages 557--563, 2006.

\bibitem[Anderson et~al.(1987)Anderson, Olkin, and
  Underhill]{doi:10.1137/0908055}
T.~Anderson, I.~Olkin, and L.~Underhill.
\newblock Generation of random orthogonal matrices.
\newblock \emph{SIAM Journal on Scientific and Statistical Computing},
  8\penalty0 (4):\penalty0 625--629, 1987.

\bibitem[Barvinok(2006)]{Barvinok}
A.~Barvinok.
\newblock Approximating orthogonal matrices by permutation matrices.
\newblock \emph{Pure and Applied Mathematics Quarterly}, 2:\penalty0 943--961,
  2006.

\bibitem[Belabbas and Wolfe(2009)]{Belabbas369}
M.-A. Belabbas and P.~J. Wolfe.
\newblock Spectral methods in machine learning and new strategies for very
  large datasets.
\newblock \emph{Proceedings of the National Academy of Sciences}, 106\penalty0
  (2):\penalty0 369--374, 2009.

\bibitem[Beylkin et~al.(1991)Beylkin, Coifman, and
  Rokhlin]{doi:10.1002/cpa.3160440202}
G.~Beylkin, R.~Coifman, and V.~Rokhlin.
\newblock Fast wavelet transforms and numerical algorithms {I}.
\newblock \emph{Communications on Pure and Applied Mathematics}, 44\penalty0
  (2):\penalty0 141--183, 1991.

\bibitem[Biloti et~al.(2013)Biloti, Matioli, and Yuan]{BILOTI201356}
R.~Biloti, L.~C. Matioli, and J.~Yuan.
\newblock A short note on a generalization of the {Givens} transformation.
\newblock \emph{Computers \& Mathematics with Applications}, 66\penalty0
  (1):\penalty0 56--61, 2013.

\bibitem[Björck and Golub(1973)]{10.2307/2005662}
Å. Björck and G.~H. Golub.
\newblock Numerical methods for computing angles between linear subspaces.
\newblock \emph{Mathematics of Computation}, 27\penalty0 (123):\penalty0
  579--594, 1973.

\bibitem[Borel(1906)]{Borel}
E.~Borel.
\newblock \emph{Introduction geometrique a quelques theories physiques}.
\newblock Gauthier-Villars, 1906.

\bibitem[{Bracewell}(1984)]{1457236}
R.~N. {Bracewell}.
\newblock The fast {Hartley} transform.
\newblock \emph{Proceedings of the IEEE}, 72\penalty0 (8):\penalty0 1010--1018,
  1984.

\bibitem[Brent and Luk(1985)]{ManySweeps}
R.~P. Brent and F.~T. Luk.
\newblock The solution of singular value and symmetric eigenvalue problems on
  multiprocessor arrays.
\newblock \emph{SIAM J. Sci. Stat. Comput.}, 6:\penalty0 69--84, 1985.

\bibitem[{Cao} et~al.(2011){Cao}, {Bachega}, and {Bouman}]{5560826}
G.~{Cao}, L.~R. {Bachega}, and C.~A. {Bouman}.
\newblock The sparse matrix transform for covariance estimation and analysis of
  high dimensional signals.
\newblock \emph{IEEE Transactions on Image Processing}, 20\penalty0
  (3):\penalty0 625--640, 2011.

\bibitem[Chen and Zeng(2012)]{DCTandKLT}
H.~Chen and B.~Zeng.
\newblock New transforms tightly bounded by {DCT} and {KLT}.
\newblock \emph{IEEE Signal Processing Letters}, 19\penalty0 (6):\penalty0
  344--347, 2012.

\bibitem[Collins and Male(2011)]{CollinsMale2011}
B.~Collins and C.~Male.
\newblock The strong asymptotic freeness of {Haar} and deterministic matrices.
\newblock \emph{Annales Scientifiques de l'Ecole Normale Superieure}, 47, 2011.
\newblock \doi{10.24033/asens.2211}.

\bibitem[Daubechies et~al.(2010)Daubechies, DeVore, Fornasier, and
  Güntürk]{doi:10.1002/cpa.20303}
I.~Daubechies, R.~DeVore, M.~Fornasier, and C.~S. Güntürk.
\newblock Iteratively reweighted least squares minimization for sparse
  recovery.
\newblock \emph{Communications on Pure and Applied Mathematics}, 63\penalty0
  (1):\penalty0 1--38, 2010.

\bibitem[Drineas et~al.(2006)Drineas, Kannan, and Mahoney]{FastMCSVD}
P.~Drineas, R.~Kannan, and M.~Mahoney.
\newblock Fast {M}onte {C}arlo algorithms for matrices {II}: Computing a
  low-rank approximation to a matrix.
\newblock \emph{SIAM Journal on Computing}, 36\penalty0 (1):\penalty0 158--183,
  2006.

\bibitem[Drmac and Veselic(2008)]{doi:10.1137/05063920X}
Z.~Drmac and K.~Veselic.
\newblock New fast and accurate {Jacobi} {SVD} algorithm. {II}.
\newblock \emph{SIAM Journal on Matrix Analysis and Applications}, 29\penalty0
  (4):\penalty0 1343--1362, 2008.

\bibitem[Fino and Algazi(1976)]{1674569}
B.~J. Fino and V.~R. Algazi.
\newblock Unified matrix treatment of the fast {Walsh}-{Hadamard} transform.
\newblock \emph{IEEE Transactions on Computers}, C-25\penalty0 (11):\penalty0
  1142--1146, 1976.

\bibitem[Freksen and Larsen(2017)]{ToeplitzJL}
C.~B. Freksen and K.~G. Larsen.
\newblock On using {Toeplitz} and circulant matrices for
  {J}ohnson-{L}indenstrauss transforms.
\newblock In \emph{28th International Symposium on Algorithms and Computation},
  pages 32:1--32:12, 2017.

\bibitem[Frerix and Bruna(2019)]{EffectiveGivens}
T.~Frerix and J.~Bruna.
\newblock Approximating orthogonal matrices with effective {Givens}
  factorization.
\newblock \emph{arXiv:1905.05796}, 2019.

\bibitem[Golub and Van~Loan(1996)]{Golub1996}
G.~H. Golub and C.~F. Van~Loan.
\newblock \emph{Matrix Computations}.
\newblock Johns Hopkins University Press, 1996.

\bibitem[Greenewald and Hero(2014)]{KroneckerPCA}
K.~H. Greenewald and A.~O. Hero.
\newblock Kronecker {PCA} based spatio-temporal modeling of video for dismount
  classification.
\newblock \emph{Proceedings of SPIE - The International Society for Optical
  Engineering}, 9093, 2014.

\bibitem[Halko et~al.(2011)Halko, Martinsson, and
  Tropp]{FindingStructureWithRandomness}
N.~Halko, P.~Martinsson, and J.~Tropp.
\newblock Finding structure with randomness: Probabilistic algorithms for
  constructing approximate matrix decompositions.
\newblock \emph{SIAM Review}, 53\penalty0 (2):\penalty0 217--288, 2011.

\bibitem[Indyk et~al.(2014)Indyk, Kapralov, and
  Price]{Indyk:2014:SSF:2634074.2634110}
P.~Indyk, M.~Kapralov, and E.~Price.
\newblock ({Nearly}) sample-optimal sparse {F}ourier transform.
\newblock In \emph{Proceedings of the 25th Annual ACM-SIAM Symposium on
  Discrete Algorithms}, pages 480--499, 2014.

\bibitem[Jain and Haupt(2017)]{CirculantDRoperator}
S.~Jain and J.~Haupt.
\newblock Convolutional approximations to linear dimensionality reduction
  operators.
\newblock In \emph{IEEE International Conference on Acoustics, Speech and
  Signal Processing}, pages 5885--5889, 2017.

\bibitem[Johansson(1997)]{HowToRandomUnitary}
K.~Johansson.
\newblock On random matrices from the compact classical groups.
\newblock \emph{Annals of Mathematics}, 145\penalty0 (3):\penalty0 519--545,
  1997.

\bibitem[Kane and Nelson(2014)]{SparserJL}
D.~M. Kane and J.~Nelson.
\newblock Sparser {J}ohnson-{L}indenstrauss transforms.
\newblock \emph{J. ACM}, 61\penalty0 (1):\penalty0 4:1--4:23, 2014.

\bibitem[Karnik et~al.(2019)Karnik, Zhu, Wakin, Romberg, and
  Davenport]{KARNIK2019624}
S.~Karnik, Z.~Zhu, M.~B. Wakin, J.~Romberg, and M.~A. Davenport.
\newblock The fast {Slepian} transform.
\newblock \emph{Applied and Computational Harmonic Analysis}, 46\penalty0
  (3):\penalty0 624 -- 652, 2019.

\bibitem[Kempen(1966)]{JacobiIsLInearEarlyOn}
H.~P. Kempen.
\newblock On quadratic convergence of the special cyclic {Jacobi} method.
\newblock \emph{Numer. Math.}, 9:\penalty0 19--22, 1966.

\bibitem[Knyazev and Argentati(2002)]{doi:10.1137/S1064827500377332}
A.~V. Knyazev and M.~E. Argentati.
\newblock Principal angles between subspaces in an {A}-based scalar product:
  algorithms and perturbation estimates.
\newblock \emph{SIAM Journal on Scientific Computing}, 23\penalty0
  (6):\penalty0 2008--2040, 2002.

\bibitem[Kondor et~al.(2014)Kondor, Teneva, and Garg]{Kondor2014MMF}
R.~Kondor, N.~Teneva, and V.~K. Garg.
\newblock Multiresolution matrix factorization.
\newblock In \emph{Proceedings of the 31st International Conference on Machine
  Learning}, pages II--1620--II--1628, 2014.

\bibitem[Le~Magoarou et~al.(2018)Le~Magoarou, Gribonval, and
  Tremblay]{lemagoarou:hal-01416110}
L.~Le~Magoarou, R.~Gribonval, and N.~Tremblay.
\newblock {Approximate fast graph Fourier transforms via multi-layer sparse
  approximations}.
\newblock \emph{{IEEE Transactions on Signal and Information Processing over
  Networks}}, 4\penalty0 (2):\penalty0 407--420, 2018.

\bibitem[Lee et~al.(2008)Lee, Nadler, and Wasserman]{Treelets}
A.~B. Lee, B.~Nadler, and L.~Wasserman.
\newblock Treelets - an adaptive multi-scale basis for sparse unordered data.
\newblock \emph{Annals of Applied Statistics}, 2\penalty0 (2):\penalty0
  435--471, 2008.

\bibitem[{Makhoul}(1980)]{1163351}
J.~{Makhoul}.
\newblock A fast cosine transform in one and two dimensions.
\newblock \emph{IEEE Transactions on Acoustics, Speech, and Signal Processing},
  28\penalty0 (1):\penalty0 27--34, 1980.

\bibitem[Mathieu and LeCun(2014)]{LeCunFastApproximations}
M.~Mathieu and Y.~LeCun.
\newblock Fast approximation of rotations and {Hessians} matrices.
\newblock \emph{arXiv:1404.7195}, 2014.

\bibitem[Merchant et~al.(2018)Merchant, Vatwani, Chattopadhyay, Raha, Nandy,
  Narayan, and Leupers]{GGR}
F.~Merchant, T.~Vatwani, A.~Chattopadhyay, S.~Raha, S.~Nandy, R.~Narayan, and
  R.~Leupers.
\newblock Efficient realization of {Givens} rotation through
  algorithm-architecture co-design for acceleration of {QR} factorization.
\newblock \emph{arXiv:1803.05320}, 2018.

\bibitem[Rath(1982)]{Rath1982}
W.~Rath.
\newblock Fast {Givens} rotations for orthogonal similarity transformations.
\newblock \emph{Numerische Mathematik}, 40\penalty0 (1):\penalty0 47--56, 1982.

\bibitem[Rusu and Thompson(2017)]{FastSparsifyingTransforms}
C.~Rusu and J.~Thompson.
\newblock Learning fast sparsifying transforms.
\newblock \emph{IEEE Transactions on Signal Processing}, 65\penalty0
  (16):\penalty0 4367--4378, 2017.

\bibitem[Schonemann(1966)]{Proc}
P.~Schonemann.
\newblock A generalized solution of the orthogonal {Procrustes} problem.
\newblock \emph{Psychometrika}, 31\penalty0 (1):\penalty0 1--10, 1966.

\bibitem[Shalit and Chechik(2014)]{Shalit14}
U.~Shalit and G.~Chechik.
\newblock Coordinate-descent for learning orthogonal matrices through givens
  rotations.
\newblock In \emph{Proceedings of the 31st International Conference on Machine
  Learning}, pages I--548--I--556, 2014.

\bibitem[Simon(2007)]{SIMON2007120}
B.~Simon.
\newblock {CMV} matrices: Five years after.
\newblock \emph{Journal of Computational and Applied Mathematics}, 208\penalty0
  (1):\penalty0 120--154, 2007.

\bibitem[Sirovich and Kirby(1987)]{Sirovich:87}
L.~Sirovich and M.~Kirby.
\newblock Low-dimensional procedure for the characterization of human faces.
\newblock \emph{J. Opt. Soc. Am. A}, 4\penalty0 (3):\penalty0 519--524, 1987.

\bibitem[Spengler et~al.(2010)Spengler, Huber, and Hiesmayr]{Spengler_2010}
C.~Spengler, M.~Huber, and B.~C. Hiesmayr.
\newblock A composite parameterization of unitary groups, density matrices and
  subspaces.
\newblock \emph{Journal of Physics A: Mathematical and Theoretical},
  43\penalty0 (38):\penalty0 385306, 2010.

\bibitem[Stewart(2019)]{Stewart2019}
K.~Stewart.
\newblock Total variation approximation of random orthogonal matrices by
  {Gaussian} matrices.
\newblock \emph{Journal of Theoretical Probability}, 2019.
\newblock ISSN 1572-9230.
\newblock \doi{10.1007/s10959-019-00900-5}.

\bibitem[Strang(2010)]{Strang12413}
G.~Strang.
\newblock Fast transforms: banded matrices with banded inverses.
\newblock \emph{Proceedings of the National Academy of Sciences}, 107\penalty0
  (28):\penalty0 12413--12416, 2010.

\bibitem[Tilma and Sudarshan(2002)]{Tilma:2002ke}
T.~E. Tilma and G.~Sudarshan.
\newblock {Generalized {Euler} angle parametrization for {SU(N)}}.
\newblock \emph{J. Phys.}, A35:\penalty0 10467--10501, 2002.

\bibitem[Tropp(2015)]{10.1561/2200000048}
J.~A. Tropp.
\newblock An introduction to matrix concentration inequalities.
\newblock \emph{Found. Trends Mach. Learn.}, 8\penalty0 (1--2):\penalty0
  1--230, 2015.

\bibitem[Van~Loan(1992)]{doi:10.1137/1.9781611970999}
C.~Van~Loan.
\newblock \emph{Computational Frameworks for the Fast {Fourier} Transform}.
\newblock Society for Industrial and Applied Mathematics, 1992.

\bibitem[Wang et~al.(2014)Wang, Lu, and Liu]{MINIMAXOPTIMALSPCA}
Z.~Wang, H.~Lu, and H.~Liu.
\newblock Tighten after relax: Minimax-optimal sparse {PCA} in polynomial time.
\newblock \emph{Advances in neural information processing systems}, pages
  3383--3391, 2014.

\bibitem[Wickerhauser(1994)]{WaveletsForPCA}
M.~V. Wickerhauser.
\newblock Two fast approximate wavelet algorithms for image processing,
  classification, and recognition.
\newblock \emph{Optical Engineering}, 33:\penalty0 33 -- 33 -- 11, 1994.

\bibitem[Wisdom et~al.(2016)Wisdom, Powers, Hershey, Roux, and Atlas]{URNN}
S.~Wisdom, T.~Powers, J.~R. Hershey, J.~L. Roux, and L.~E. Atlas.
\newblock Full-capacity unitary recurrent neural networks.
\newblock In \emph{Advances in Neural Information Processing Systems 29}, pages
  4880--4888, 2016.

\bibitem[{Zhang} et~al.(2014){Zhang}, {Liu}, {Dong}, and {Wang}]{SparseFaces}
H.~{Zhang}, W.~{Liu}, L.~{Dong}, and Y.~{Wang}.
\newblock Sparse eigenfaces analysis for recognition.
\newblock In \emph{12th International Conference on Signal Processing}, pages
  887--890, 2014.

\bibitem[Zou et~al.(2006)Zou, Hastie, and Tibshirani]{SparsePCA}
H.~Zou, T.~Hastie, and R.~Tibshirani.
\newblock Sparse principal component analysis.
\newblock \emph{Journal of Computational and Graphical Statistics}, 15\penalty0
  (2):\penalty0 265--286, 2006.

\end{thebibliography}
}

\newpage

\section*{Supplementary materials}

\begin{figure}[!h]
		\centering
		\includegraphics[trim = 15 5 5 14, clip, width=0.3\textwidth]{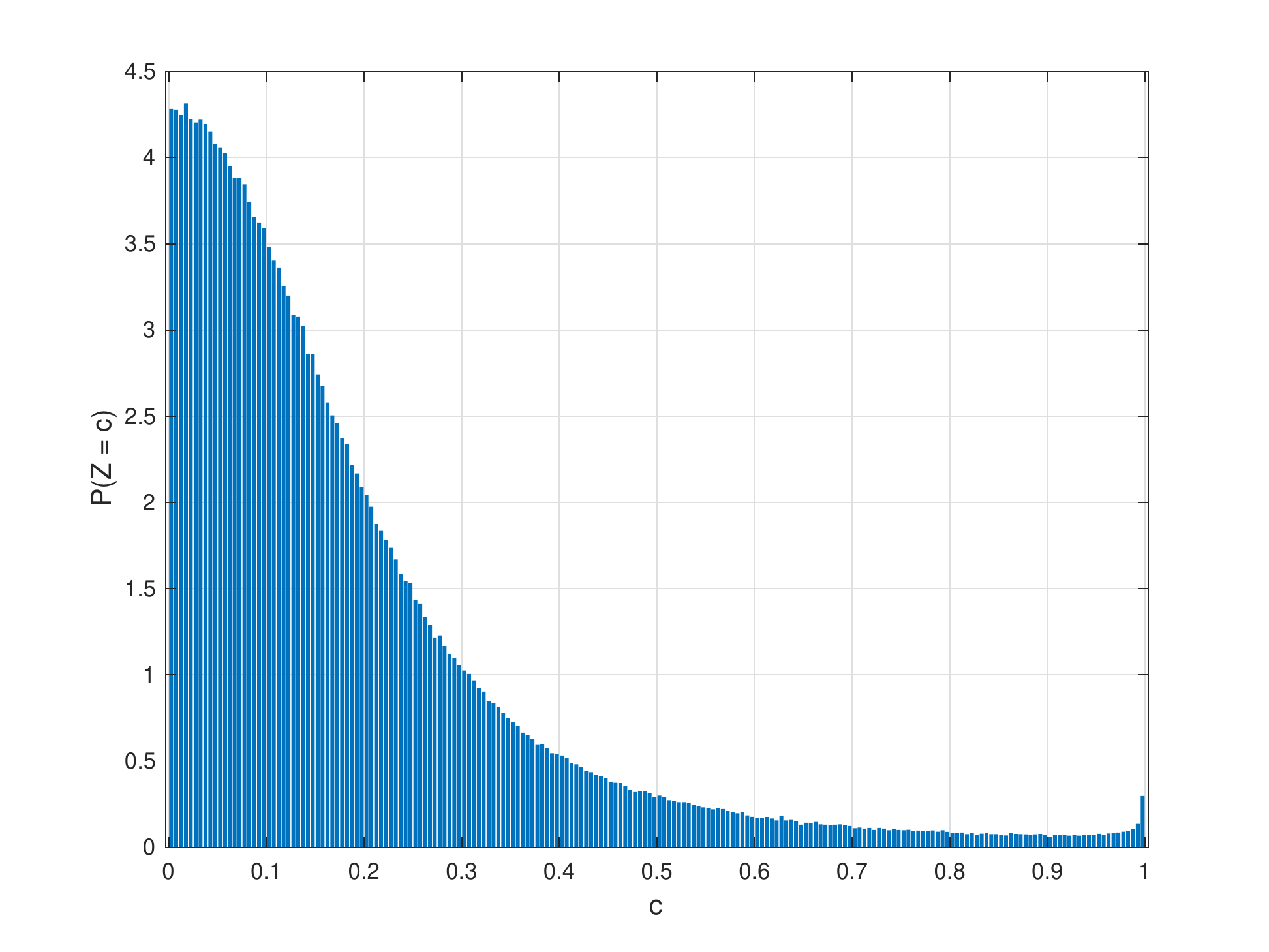}
		\includegraphics[trim = 10 5 5 14, clip, width=0.3\textwidth]{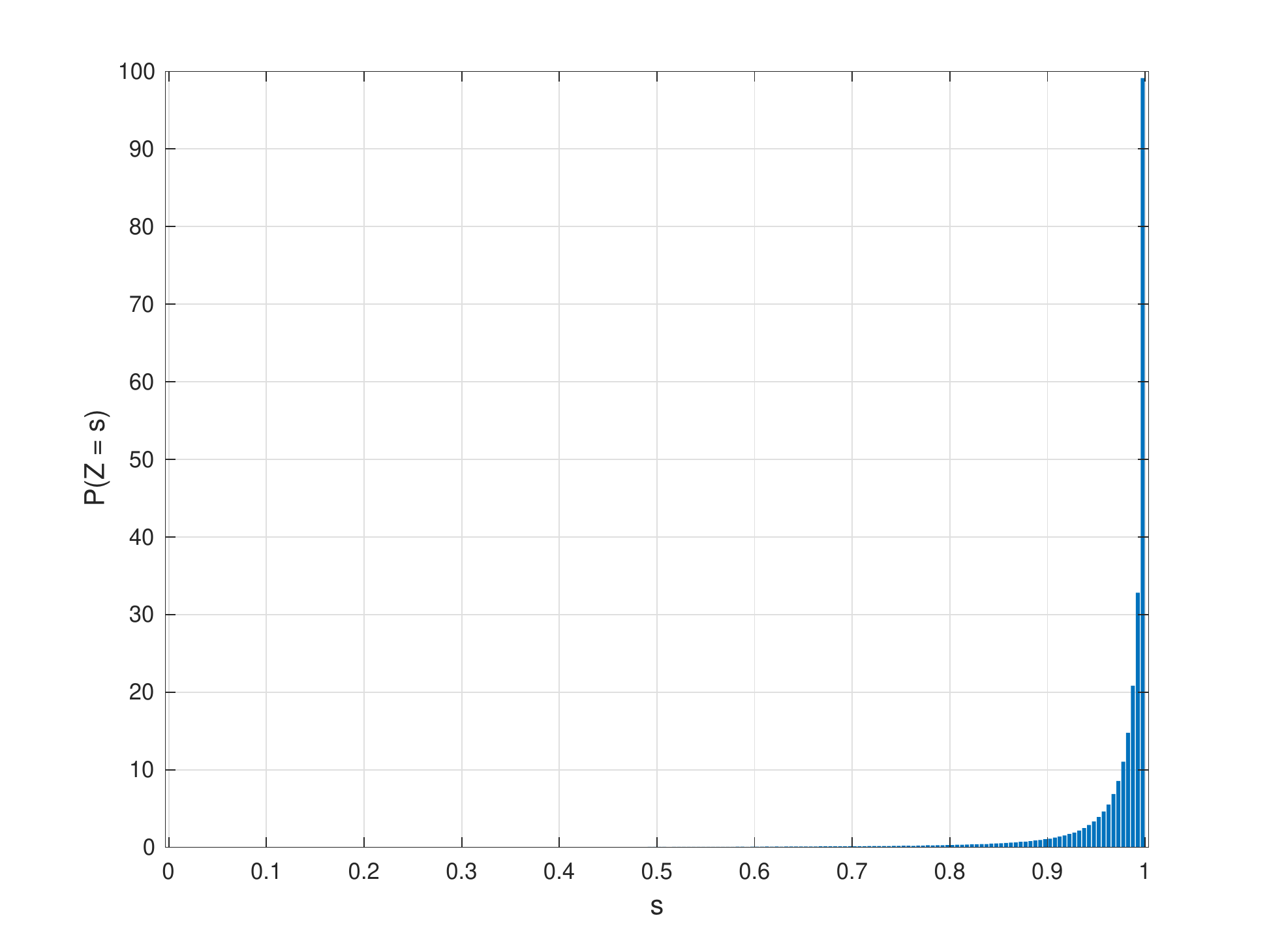}
		\caption{Experimental distribution on the entries (cosine and sine on the left and right, respectively) of the Givens rotations $\mathbf{G}_{ij}$ from \eqref{eq:therandom}. We observe numerically that $\cos \theta$ follows closely an exponential distribution $\lambda \exp (-\lambda c)$ with $\lambda \approx \sqrt{d}/2$ on the interval $[0,1]$.}
		\label{fig:moreremark3}
\end{figure}

\noindent \textbf{More on Remark 3.} We say that the orthonormal group has $O(d^2)$ degrees of freedom as it was established \cite{doi:10.1137/0908055} that any orthonormal matrix $\mathbf{U}$ can be factored into a product as
\begin{equation}
	\mathbf{U}  = \left( \prod_{i=1}^d \prod_{j=i+1}^d \mathbf{G}_{ij}(\theta_{ij})  \right) \mathbf{D},
	\label{eq:therandom}
\end{equation}
where $\mathbf{D}$ is a diagonal matrix with entries only in $\{ \pm 1 \}$ and the $\mathbf{G}_{ij}(\theta_{ij})$ are Givens rotations, with angles $\theta_{ij}$, i.e., we have $c = \cos \theta_{ij}$ and $s = \sin \theta_{ij}$ in \eqref{eq:theG}. To generate a random orthonormal $\mathbf{U}$ we therefore need to generate random $\mathbf{D}$ (which are $\{ \pm 1 \}$ with equal probability) and $\frac{d(d-1)}{2}$ random angles $0 \leq \theta_{ij} \leq \pi/2$. These angles are mutually independent and it is known that their joint density function is a random variable
\begin{equation}
	Z \propto \left( \prod_{k=2}^d \cos^{k-2} \theta_{1k}  \right) \left( \prod_{k=3}^d \cos^{k-3} \theta_{2k}  \right) \left( \prod_{k=d}^d \cos^{k-d} \theta_{(d-1)k}  \right).
	\label{eq:jointdensity}
\end{equation}
We define the beta random variable $\sqrt{y} = \cos \theta$ (and $\pm \sqrt{1-y} = \sin \theta$) with density
\begin{equation}
	f(y, \alpha, \beta) = \frac{ y^{\alpha-1} (1-y)^{\beta-1} }{B(\alpha, \beta)}, \ B(\alpha, \beta) = \frac{\Gamma(\alpha) \Gamma(\beta)}{\Gamma(\alpha+\beta)} \text{ with } y \in [0,1].
\end{equation}
Because we are interested in the computational complexity of a random orthonormal matrix we focus on the following three special cases: i) do nothing: $\mathbf{\tilde{G}}_{ij} = \begin{bmatrix} \pm 1 & 0 \\ 0 & \pm 1\end{bmatrix}$; ii) permute coordinates: $\mathbf{\tilde{G}}_{ij} = \begin{bmatrix} 0 & \pm 1 \\ \pm 1 & 0 \end{bmatrix}$; and iii) $\mathbf{\tilde{G}}_{ij} \in \left\{ \frac{1}{\sqrt{2}} \begin{bmatrix} 1 & 1 \\ -1 & 1 \end{bmatrix},  \frac{1}{\sqrt{2}} \begin{bmatrix} 1 & -1 \\ 1 & 1 \end{bmatrix} \right\}$. The first two cases perform no operations while the last performs only 4 (as compared to 6 for a general rotation). Unfortunately, the joint density in \eqref{eq:jointdensity} does not seem to have a simlpe closed-form expression. As such we show in Figure \ref{fig:moreremark3} numerical results of the probability distribution of $\cos \theta$ over all angles $\theta_{ij}$, which we observe numerically that approaches an exponential distribution $\lambda \exp(-\lambda c )$. Concentration around the three special cases does not occur and therefore a random orthonormal matrix will generally have computational complexity $O(d^2)$. For example, if we discretized the continuum of $c = \cos \theta$ then the probability that a random $\mathbf{U}$ is an approximate permutation matrix and therefore basically exhibits no numerical complexity is $(1-\exp(-\lambda \epsilon ) )^{\frac{d(d-1)}{2}}$ for $0 < \epsilon \ll 1$, i.e., the probability that all $\frac{d(d-1)}{2}$ rotations have $\cos \theta \leq \epsilon$.$\hfill \square$

\noindent \textbf{Proof of Theorem 1.} For simplicity of exposition we will drop the sub-index $k$ herein and therefore \eqref{eq:thestart} develops to the following
\begin{equation}
\| \mathbf{L} - \mathbf{G}_{i j} \mathbf{N} \|_F^2 \! = \!  \| \mathbf{L} \|_F^2 \! + \! \| \mathbf{G}_{ij} \mathbf{N} \|_F^2 \! - \! 2\text{tr}(\mathbf{N}^T \mathbf{G}_{i j}^T \mathbf{L})
= \| \mathbf{\sigma}_p \|_2^2 + \! \| \mathbf{\bar{\sigma}}_p \|_2^2 - 2\text{tr}( \mathbf{G}_{i j}^T \mathbf{L} \mathbf{N}^T),
%& = \! \! \| \mathbf{\sigma}_p \|_2^2 + \! \| \mathbf{\bar{\sigma}}_p \|_2^2 \! - \! 2\text{tr}(\mathbf{Z}) - 2(\| \mathbf{Z}_{\{i, j\}} \|_* \! \! - \! \text{tr}(\mathbf{Z}_{\{i, j\}} \!) \\
%& = \! \! \| \mathbf{\sigma}_p \|_2^2 + \! \| \mathbf{\bar{\sigma}}_p \|_2^2 - 2\text{tr}(\mathbf{Z}) - 2C_{i j},
\label{eq:givensapproxforortho}
\end{equation}
where the Frobenius norms reduces to the $\ell_2$ norms of the spectra and we have used the circular permutation property of the trace. It is convenient to denote $\mathbf{Z} = \mathbf{L} \mathbf{N}^T$ and the $2 \times 2$ matrix $\mathbf{Z}_{\{i,j\}} = \begin{bmatrix} Z_{ii} & Z_{ij}  \\ Z_{ji} & Z_{jj} \end{bmatrix} $. Given that $\mathbf{G}_{ij}$ performs operations only on rows $i$ and $j$, the trace is
\begin{equation}
\text{tr}(\mathbf{G}_{ij}^T \mathbf{L}\mathbf{N}^T) = \sum_{k=1, k \notin \{i,j\}}^d Z_{kk} + \text{tr}(\mathbf{\tilde{G}}_{ij}^T \mathbf{Z}_{ \{ i,j \} })
% \text{tr}(\mathbf{Z}) - Z_{ii} - Z_{jj} + \text{tr}(\mathbf{\tilde{G}}_{ij}^T \mathbf{Z}_{ \{ i,j \} }) \\
= \text{tr}(\mathbf{Z}) + \text{tr}(\mathbf{\tilde{G}}_{ij}^T \mathbf{Z}_{ \{ i,j \} }) - \text{tr}(\mathbf{Z}_{ \{i,j\} } ).
\label{eq:thetraceteerm}
\end{equation}
To minimize the quantity in \eqref{eq:givensapproxforortho} we have to maximize \eqref{eq:thetraceteerm} which is known as a Procrustes problem \cite{Proc} whose solution is given by the polar decomposition of $\mathbf{Z}_{ \{i,j\} }$ detailed in \cite{Golub1996}[Chapter~9.4.3]. Therefore, we set the optimal transformation to
\begin{equation}
\mathbf{\tilde{G}}_{i j}^\star = \mathbf{V}_1 \mathbf{V}_2^T,\ \mathbf{Z}_{\{i,j\}} = \mathbf{V}_1 \mathbf{S V}_2^T,
\label{eq:besttransformation}
\end{equation}
where use the SVD of $\mathbf{Z}_{ \{i,j\} }$ ($\mathbf{S} = \text{diag}(s_1, s_2)$ are the singular values). With this choice, we have
\begin{equation}
\text{max} \  \text{tr}(\mathbf{G}_{ij}^T \mathbf{L} \mathbf{N}^T \! ) \! %= \!  \text{tr}(\mathbf{Z}) \! \! +\! \text{tr}(\!\mathbf{\tilde{G}}_{i j}^{\star T} \mathbf{Z}_{ \{ i,j \} }\! ) \! - \! \text{tr}(\!\mathbf{Z}_{ \{ i,j \} } \! )
%= \text{tr}(\mathbf{Z}) + \text{tr}(\mathbf{V}_2 \mathbf{V}_1^T \mathbf{V}_1 \mathbf{SV}_2^T) - \text{tr}(\mathbf{Z}_{ \{ i,j \} })
= \text{tr}(\mathbf{Z}) + \text{tr}(\mathbf{S}) - \text{tr}(\mathbf{Z}_{ \{ i,j \} })
= \text{tr}(\mathbf{Z}) + \| \mathbf{Z}_{ \{i,j\} } \|_* - \text{tr}(\mathbf{Z}_{ \{ i,j \} }) 
=  \text{tr}(\mathbf{Z}) +C_{ij}.
\label{eq:somelongdevelopment}
\end{equation}
We denote the nuclear norm $\| \mathbf{Z}_{\{ i,j\}} \|_*$, i.e., the sum of the singular values $s_1$ and $s_2$ and we define
\begin{equation}
C_{ij} = \| \mathbf{Z}_{ \{i,j\} } \|_* - \text{tr}(\mathbf{Z}_{ \{ i,j \} }).
\label{eq:theCija}
\end{equation}
Intuitively, the results \eqref{eq:givensapproxforortho}, \eqref{eq:thetraceteerm} follows after observing that: 1) the $\mathbf{G}_{ij}$ can be viewed as a perturbed identity matrix; 2) if $\mathbf{G}_{ij}$ is exactly $\mathbf{I}_{d \times d}$ then $\| \mathbf{L} - \mathbf{N} \|_F^2 = \| \mathbf{\sigma}_p \|_2^2 + \| \mathbf{\bar{\sigma}}_p \|_2^2 - 2\text{tr}(\mathbf{Z})$ while if $\mathbf{G}_{ij}$ is the optimal orthonormal transformation $\mathbf{Q}$ that minimizes \eqref{eq:givensapproxforortho} given by the Procrustes solution \citep{Proc} then we have $\| \mathbf{L} - \mathbf{Q} \mathbf{N} \|_F^2 = \| \mathbf{\sigma}_p \|_2^2 + \| \mathbf{\bar{\sigma}}_p \|_2^2 - 2\| \mathbf{Z} \|_*$, where the last term is the nuclear norm of $\mathbf{Z}$; 3) therefore, we actually apply the identity transformation on all coordinates, i.e., the $\text{tr}(\mathbf{Z})$ term, while for the two chosen coordinates we apply the best (in the sense of reducing the error) orthogonal transformation whose contribution is the nuclear norm term $\| \mathbf{Z}_{\{ i,j \} } \|_*$ and then correct for the trace term that was wrongly added initially in $\text{tr}(\mathbf{Z})$, by subtracting $\text{tr}(\mathbf{Z}_{ \{i,j \}})$.\\
There are $d(d-1)/2$ quantities $C_{ij}$ but they can be computed efficiently by noting that the singular values of $\mathbf{Z}_{ \{i,j\} }$ are $s_{1,2} \! = \!  \sqrt{ \! \frac{1}{2} \! \! \left( \!   \| \mathbf{Z}_{\{ i,j\}} \|_F^2 \!  \pm  \! \sqrt{\! \|\mathbf{Z}_{\{ i,j\}}\|_F^4 \!- \! 4 \det(\mathbf{Z}_{\{ i,j\}})^2} \right) }$. Observe that both singular values are of the form $\sigma_{1,2} = \sqrt{A \pm \sqrt{B}}$ which can be written as $\sqrt{X} \pm \sqrt{Y}$ where $X = \frac{A + \sqrt{A^2-B}}{2}$ and $Y = A - X$. Written like this we can see that $\| \mathbf{Z}_{ \{i,j \} } \|_* = \sigma_1 + \sigma_2 = 2X = \sqrt{ \| \mathbf{Z}_{ \{i,j\} }  \|_F^2 + 2 | \det( \mathbf{Z}_{ \{i,j\} } ) | }$. Depending on the sign of the determinant we have either $\| \mathbf{Z}_{ \{i,j \} } \|_* = \sqrt{ (Z_{ii} + Z_{jj})^2 + (Z_{ij} - Z_{ji})^2 }$ or $\| \mathbf{Z}_{ \{i,j \} } \|_* = \sqrt{ (Z_{ii} - Z_{jj})^2 + (Z_{ij} + Z_{ji})^2 }$, respectively. These give the final formulas for $C_{ij}$ from \eqref{eq:theCijcomputed}.
$\hfill \blacksquare$

\noindent \textbf{Proof of Theorem 2.} Assume  $d$ is even and  partition the set $\{1, \dots, d\}$ into $d/2$ pairs of indices $(i_k, j_k)$. We therefore have
\begin{equation}
\sum_{k=1}^{d/2} C_{i_k j_k} = \sum_{k=1}^{d/2} ( \| \mathbf{U}_{ \{i_k, j_k\} } \|_*  - U_{i_k i_k} - U_{j_k j_k} )
= \sum_{k=1}^{d/2} \| \mathbf{U}_{ \{i_k, j_k\} } \|_* - \text{tr}(\mathbf{U}).
\end{equation}
As a side note, to maximization of this partitioned quantity is related to the weighted maximum matching algorithm (of maximum-cardinality matchings) on the graph with $d$ nodes and with edge weights $C_{i_k j_k}$. With this choice of indices, the objective function becomes
\begin{equation}
\left\| \mathbf{U} - \prod_{k=1}^{d/2} \mathbf{G}_{i_k j_k} \right\|_F^2 = 2d - 2\text{tr}(\mathbf{U}) - 2 \sum_{k=1}^{d/2} C_{i_k j_k} =  2d -2\sum_{k=1}^{d/2} \| \mathbf{U}_{ \{i_k j_k\} } \|_*.
\label{eq:themanyG}
\end{equation}
We use the singular value decomposition of a generic $\mathbf{U}_{ \{i, j\} } = \mathbf{V}_1 \mathbf{S V}_2^T, \mathbf{S} = \text{diag}(\mathbf{s})$, and develop:
\begin{equation}
|\text{tr}( \mathbf{U}_{ \{i, j\} })| \! \! = \! \! |\text{tr}(\mathbf{V}_1\mathbf{S V}_2^T)| \! \! =\! \! |\text{tr}(\mathbf{S V}_2^T \mathbf{V}_1)| \! \! = \! \!   |\text{tr}(\mathbf{S \Delta})|\! =\! \Biggl\lvert \sum_{t=1}^2 \! s_{t} \Delta_{tt} \Biggr\rvert \! \leq \! \Delta_{\max} \! \sum_{t=1}^2 s_t \! \! = \! \! \Delta_{\max} \| \mathbf{U}_{ \{i, j\} } \|_*,
\label{eq:traceandnuclear}
\end{equation}
where we have use the circular property of the trace and $\mathbf{\Delta} = \mathbf{V}_2^T \mathbf{V}_1$ where $\Delta_{tt}$ are its diagonal entries which obey $| \Delta_{tt} | \leq \Delta_{\max}$. We define the diagonal coherence as
\begin{equation}
	\Delta_{\max} = \max \{ | \Delta_{11} |, | \Delta_{22} | \}.
\end{equation}
With \eqref{eq:traceandnuclear}, we can state that for the $k^\text{th}$ transformation  that
\begin{equation}
\mathbb{E}[\| \mathbf{U}_{ \{i_k, j_k\} } \|_*] \geq \frac{\pi}{2} \mathbb{E}[\text{tr}(\mathbf{U}_{ \{i_k, j_k\} })],
\label{eq:nucleartracebound}
\end{equation}
where we have used the fact that $\mathbf{V}_1$ and $\mathbf{V}_2$ have the structure $\mathbf{\tilde{G}}_{i j}$ in \eqref{eq:theG} and therefore
\begin{equation*}
\mathbb{E}[\Delta_{\max}] \! \!  = \! \! \frac{1}{\pi^2} \! \!  \iint_{0}^{\pi} \! \! \! \! \! |\cos(x)\cos(y) +\sin(x)\sin(y)| dx dy \! = \! \frac{2}{\pi}.
\end{equation*}
Finally, given an orthonormal $\mathbf{U}$ of size $d \times d,\ d \geq 4$, we use \eqref{eq:themanyG} and \eqref{eq:nucleartracebound} to bound
\begin{equation}
\begin{aligned}
\mathbb{E} \left[ \! \left\| \mathbf{U} - \prod_{k=1}^{d/2} \mathbf{G}_{i_k j_k} \right\|_F^2 \! \right] & =  2d - 2 \sum_{k=1}^{d/2} \! \mathbb{E}[\| \mathbf{U}_{ \{i_k, j_k\} } \|_*]  \leq 2d - \pi \sum_{k=1}^{d/2} \mathbb{E}[\text{tr} (\mathbf{U}_{ \{i_k, j_k\} })] \\
& \leq 2d - \pi \mathbb{E}[\text{tr}(\mathbf{U})] = 2d - \pi \mathbb{E} \left[\sum_{t=1}^{d} | U_{tt} | \right] = 2d - \sqrt{2 \pi d},
\end{aligned}
\label{eq:interesting}
\end{equation}
where we have used that $\mathbb{E} \left[\sum_{t=1}^{d} | U_{tt} | \right] = \sqrt{2 \pi^{-1} d}$ because the diagonal elements of $\mathbf{U}$ can be viewed as Gaussian random variables with zero mean and standard deviation $d^{-1/2}$ (as the columns of $\mathbf{U}$ are normalized in the $\ell_2$ norm) \cite{Stewart2019} and because the $\ell_1$ norm of a standard Gaussian random vector of size $d$ is $\sqrt{2\pi^{-1}}d$.\\
When $d$ is odd, we extend the matrix with a zero column/row and ``1'' on the diagonal and the argument follows  in the same  way. In \eqref{eq:themanyG} we could use the expected value calculated in \eqref{eq:expectedvalue} but we reach a worse, lower, constant in \eqref{eq:interesting} for the $-\sqrt{d}$ term and therefore a worse overall bound.$\hfill \blacksquare$

\noindent \textbf{Proof of Theorem 3.} Given the orthonormal $\mathbf{U}$, by \cite{Golub1996}[Theorem~5.2.1], we can construct its QR factorization using a set of Givens rotations \cite{Golub1996}[Chapter~5.2.5]. After introducing zeros in the first $r$ columns of $\mathbf{U}$, by left multiplication with Givens rotations, we reach its following partial triangularization
\begin{equation}
\mathbf{J}_{i_g j_g} \dots \mathbf{J}_{i_1 j_1} \mathbf{U} = \begin{bmatrix}
\mathbf{D}'' & \mathbf{0}_{r \times (d-r)}\\
\mathbf{0}_{(d-r) \times r} & \mathbf{U}'
\end{bmatrix},
\label{eq:partialtriang}
\end{equation}
where the diagonal matrix $\mathbf{D}''$ of size $r \times r$ has entries $D''_{tt} \in \{ \pm 1 \}$ and $\mathbf{U}'$ of size $(d-r) \times (d-r)$ is orthonormal. To introduce the zeros on the $t^\text{th}$ column we need $(d-t)$ Givens rotations and therefore to bring $\mathbf{U}$ to the structure in \eqref{eq:partialtriang} we need $g = \frac{r}{2}(2d - r - 1)$ Givens rotations which we have denoted $\mathbf{J}_{i_k j_k},\ k = 1,\dots,g$. We are exploiting the fact that the triangularization of an orthogonal matrix leads to a $\pm 1$ diagonal. Then we might consider a good approximation to $\mathbf{U}$ the product $\mathbf{\bar{U}} = \mathbf{J}_{i_1 j_1}^T \dots \mathbf{J}_{i_g j_g}^T \mathbf{D}$.
where $\mathbf{D} = \begin{bmatrix} \mathbf{D}'' & \mathbf{0}_{r \times (d-r)} \\ \mathbf{0}_{(d-r) \times r} & \mathbf{D}' \end{bmatrix}$ with $D'_{tt} = \text{sgn}(U'_{tt})$ and $\mathbf{D}''$ is taken from \eqref{eq:partialtriang}. The goal of the diagonal matrix $\mathbf{D}$ is to ensure that the product in $\mathbf{\bar{U}}$ has a nonnegative diagonal. Then, given $g$ transforms  we can bound
\begin{equation}
\| \mathbf{U} - \mathbf{J}_{i_1 j_1}^T \dots \mathbf{J}_{i_g j_g}^T \mathbf{D} \|_F^2 =  2(d - \lfloor r \rfloor) - 2\text{tr}(\mathbf{D}'\mathbf{U}').
\end{equation}
If we consider $\mathbf{\tilde{G}}_{i j}$ in \eqref{eq:theG} instead of the rotations $\mathbf{J}_{i_k j_k}$ then the quantity on the right becomes an upper bound, since Givens rotations are a special case of $\mathbf{\tilde{G}}_{i j}$ -- we can always initialize the $\mathbf{G}_{i_k j_k}$ of Algorithm 1 with the $\mathbf{J}_{i_k j_k}$ defined above and the iterative procedure is guaranteed not to worsen the factorization. Therefore, the result follows after using $\mathbb{E} \left[\text{tr}(\mathbf{D}'\mathbf{U}') \right] = \mathbb{E} \left[\sum_{t=1}^{d-\lfloor r \rfloor} | U'_{tt} | \right] = \sqrt{2 (d - \lfloor r \rfloor ) \pi^{-1}}$.$\hfill \blacksquare$

\noindent \textbf{Proof of Theorem 4.} First, we introduce the off-diagonal ``norm'', i.e., the square-root of the squared sum of the off-diagonal elements of an orthonormal matrix $\mathbf{U} \in \mathbb{R}^{d \times d}$ as
\begin{equation}
\text{off}(\mathbf{U})^2  =  \sum_{t=1}^d \sum_{q=1, q \neq t}^d U_{tq}^2 = \| \mathbf{U} \|_F^2 - \sum_{t=1}^d U_{tt}^2 = d - \sum_{t=1}^d U_{tt}^2.
\label{eq:offU}
\end{equation}
Better approximations $\mathbf{\bar{U}}$ of $\mathbf{U}$ lead to lower $\text{off}(\mathbf{U} \mathbf{\bar{U}}^T)$, as $\mathbf{U} \mathbf{\bar{U}}^T$ approaches the identity. If we use this measure, we reach the following result.\\
We use the fact that $C_{ij}$ is added to the diagonal of $\mathbf{U}$, but only to $U_{ii}$ and $U_{jj}$. The quantity $f(\gamma) = (U_{ii} + \gamma C_{ij})^2 + (U_{jj} + (1- \gamma) C_{ij})^2, \gamma \in \mathbb{R},$ is minimized for $C_{ij} \neq 0$ when $\gamma_0 = \frac{1}{2} + \frac{U_{jj}-U_{ii}}{2C_{ij}}$  (and therefore $f(\gamma_0) = \frac{1}{2}\| \mathbf{U}_{ \{ i,j \} } \|_*^2 $) which leads to
\begin{equation}
\begin{aligned}
&\text{off}(\mathbf{UG}_{ij}^T )^2 = d - \sum_{t=1, t \notin \{ i, j\} }^d U_{tt}^2 - f(\gamma)
\leq d - \sum_{t=1, t \notin \{ i, j\} }^d U_{tt}^2 - f(\gamma_0) \\
= & d - \sum_{t=1}^d U_{tt}^2 + U_{ii}^2 + U_{jj}^2  - f(\gamma_0) = \text{off}(\mathbf{U})^2 + \frac{(U_{ii}  - U_{jj})^2}{2} - C_{ij}(U_{ii} \! + \! U_{jj}) - \frac{C_{ij}^2}{2} \\
= & \text{off}(\mathbf{U})^2 \! \! + \! \! \frac{(U_{ii} \! \! - \! U_{jj})^2}{2} \! -  \! \frac{\| \mathbf{U}_{ \{i,j\} } \|_*^2  \! \! - \! (U_{ii} \! + \! U_{jj})^2}{2} = \text{off}(\mathbf{U})^2 + U_{ii}^2 + U_{jj}^2  - \frac{\| \mathbf{U}_{ \{i,j\} } \|_*^2}{2}.
\end{aligned}
\label{eq:nice}
\end{equation}
We now use de explicit formulas for the singular values of a $2 \times 2$ matrix and the fact that 
\begin{equation}
\| \mathbf{U}_{ \{i,j\} } \|_*^2 = s_1^2 + s_2^2 + 2s_1 s_2 \\
=  \| \mathbf{U}_{ \{i,j\} } \|_F^2 + 2|\text{det}(\mathbf{U}_{ \{i,j\} })|,
\end{equation}
and expand this expression to get in \eqref{eq:nice} 
\begin{equation*}
\begin{aligned}
\text{if }& \text{det}(\mathbf{U}_{ \{i,j\} }) \geq 0:
&\text{off}(\mathbf{UG}_{ij}^T )^2 \leq \text{off}(\mathbf{U})^2 + \frac{(U_{ii} - U_{jj})^2-(U_{ij} - U_{ji})^2}{2}, \\
\text{if }& \text{det}(\mathbf{U}_{ \{i,j\} })< 0:
&\text{off}(\mathbf{UG}_{ij}^T )^2 \leq \text{off}(\mathbf{U})^2 + \frac{(U_{ii} + U_{jj})^2 - (U_{ij} + U_{ji})^2}{2}.
\end{aligned}
\end{equation*}
Therefore, to guarantee $\text{off}(\mathbf{UG}_{ij}^T )^2 \leq \text{off}(\mathbf{U})^2$ we need $2(U_{ii}^2 + U_{jj}^2) \leq \| \mathbf{U}_{ \{i,j\} } \|_*^2$ which is equivalent to
\begin{equation}
\begin{aligned}
\text{if }& \text{det}(\mathbf{U}_{ \{i,j\} }) \geq 0:\ (U_{ii} - U_{jj})^2 \leq (U_{ij} - U_{ji})^2, \\
\text{if }& \text{det}(\mathbf{U}_{ \{i,j\} })< 0:\ (U_{ii} + U_{jj})^2 \leq (U_{ij} + U_{ji})^2.
\end{aligned}
\end{equation}
In this paper we assume that $\mathbf{U}$ is taken randomly from the Haar measure \cite{HowToRandomUnitary} and then modified to have positive diagonal. Therefore, we have $U_{tt} \geq 0$ for all $t$ by construction, otherwise we would just consider the update $U_{tt}  \leftarrow \text{sign}(U_{tt})U_{tt}$. Moreover, as the algorithm progresses we continue to have that $\text{det}(\mathbf{U}_{ \{i,j\} }) \geq 0$ because each $\mathbf{G}_{ij}$ adds a positive amount (the $C_{ij}$ value) to the diagonal elements $U_{ii}$ and $U_{jj}$ thus converging towards $\mathbf{I}_{d \times d}$ in the sense of \eqref{eq:towardsidentity}.\\
In a similar way we can construct a lower bound. Assuming w.l.o.g. that $U_{ii} \geq U_{jj} \geq 0$, the quantity $f(\gamma) = (U_{ii} + \gamma C_{ij})^2 + (U_{jj} + (1- \gamma) C_{ij})^2$ is maximized when $\gamma_0 = 1$ and therefore $f(\gamma_0) = (U_{ii} + C_{ij})^2 + U_{jj}^2 $ which, similarly to \eqref{eq:nice}, leads to
\begin{equation}
\text{off}(\mathbf{UG}_{ij}^T)^2 \geq \text{off}(\mathbf{U})^2 - C_{ij}(2U_{ii} + C_{ij}).\hfill \blacksquare
\end{equation}

\begin{figure}[!tbp]
	\centering
		\centering
		\includegraphics[trim = 15 0 30 15, clip, width=0.32\textwidth]{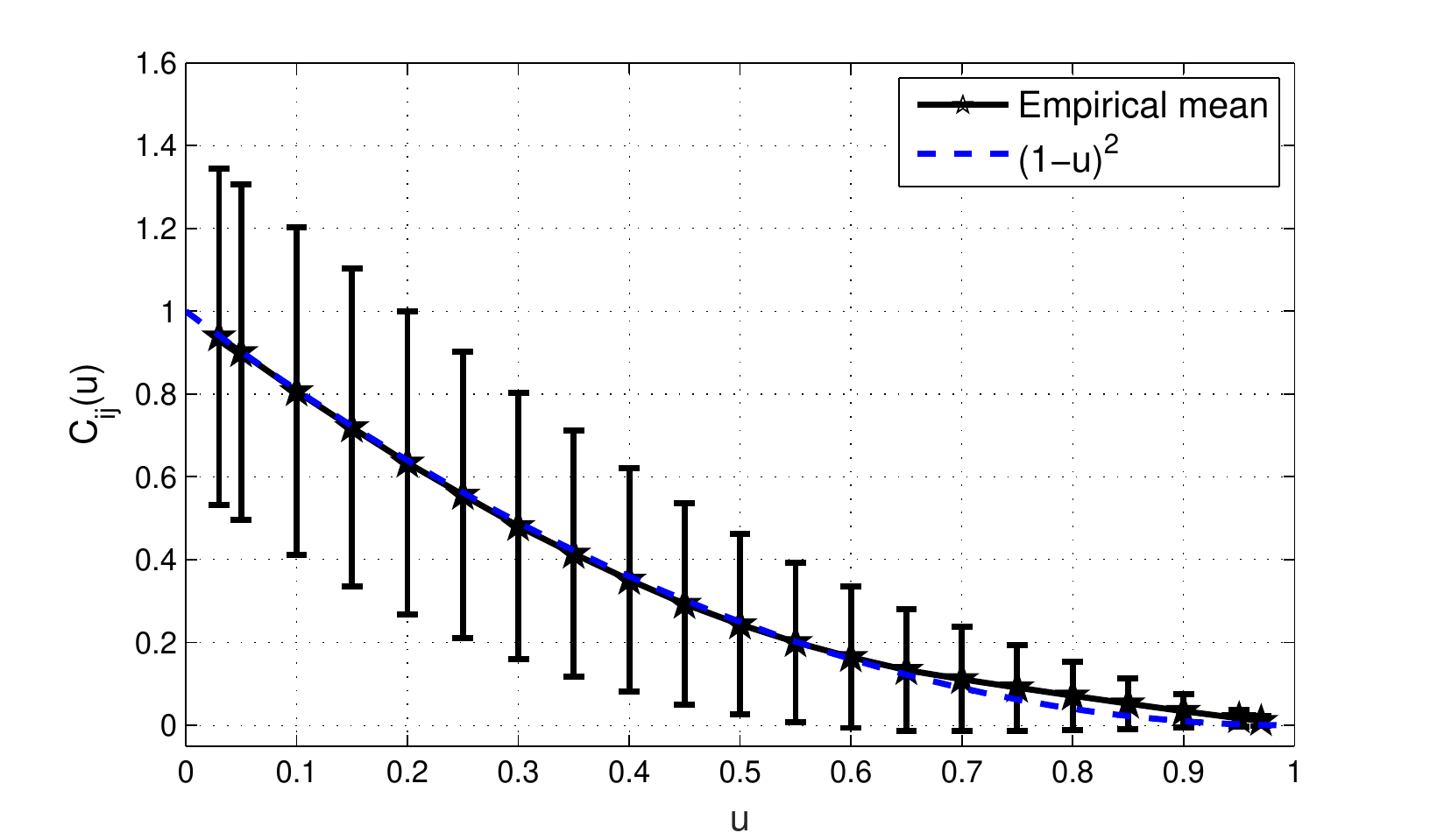}
		\includegraphics[trim = 15 5 30 20, clip, width=0.32\textwidth]{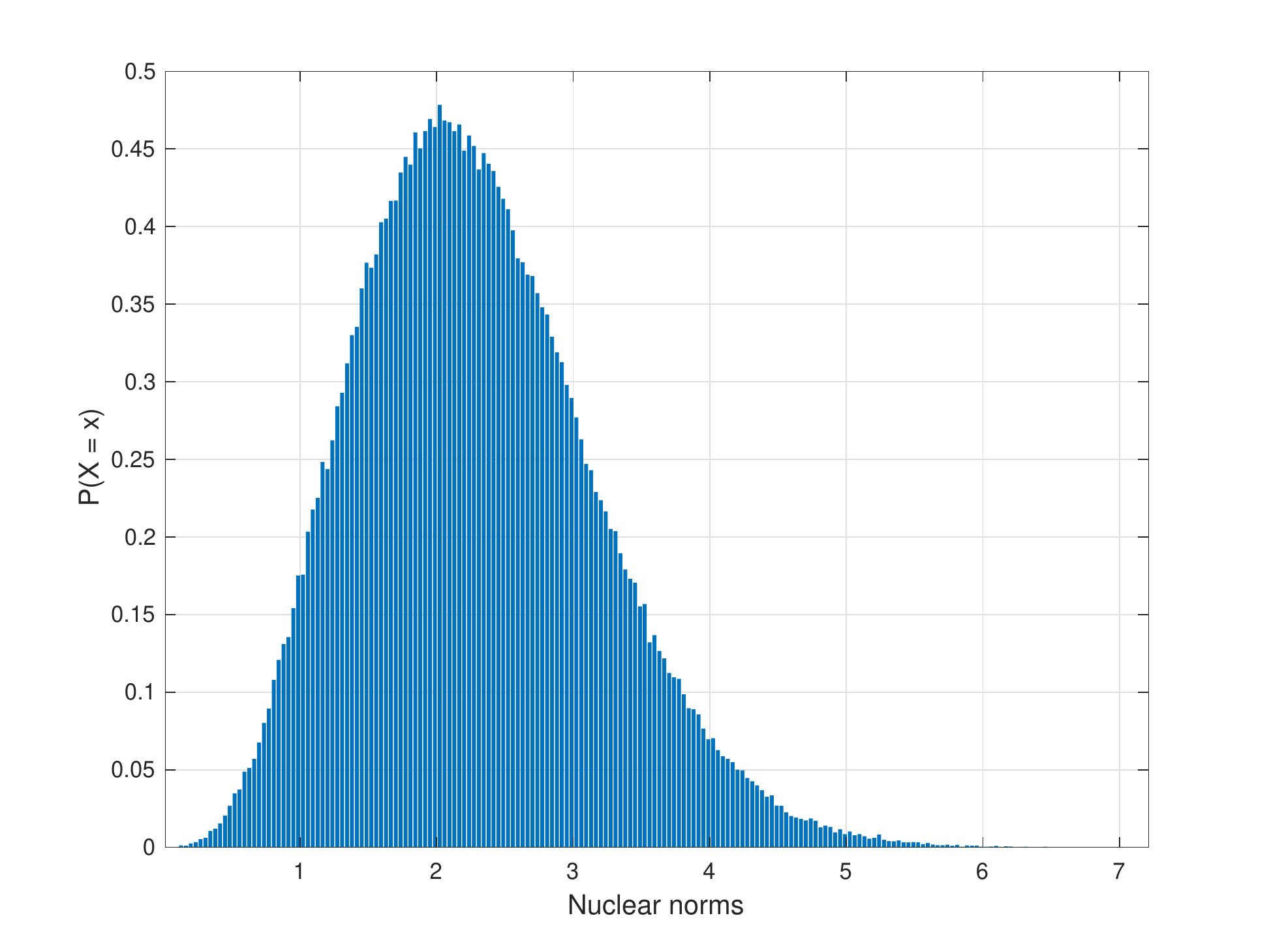}
		\includegraphics[trim = 15 5 30 20, clip, width=0.32\textwidth]{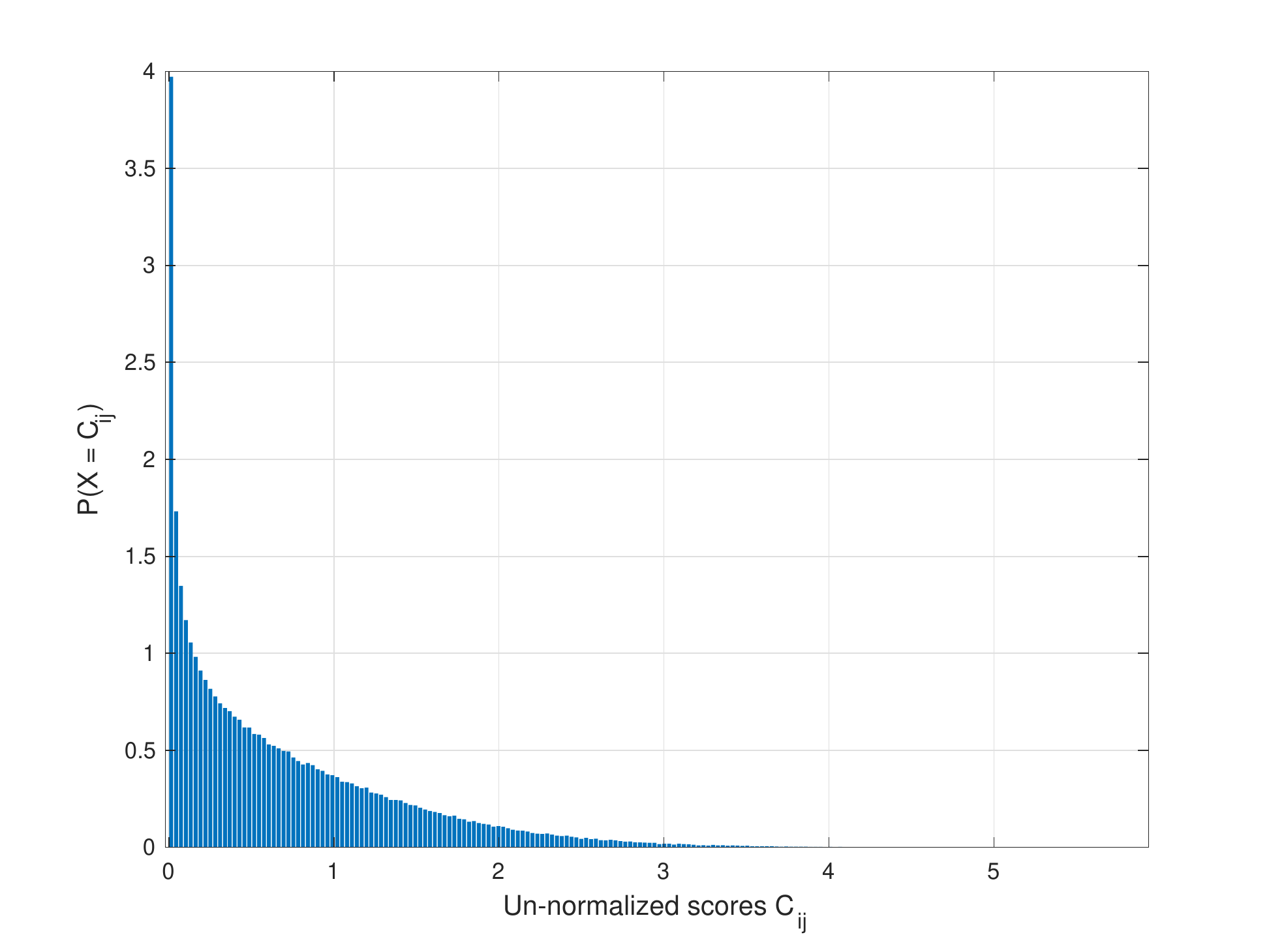}
		\caption{Left: empirical mean and standard deviation verification of Remark 5; Middle: empirical pdf for the nuclear norms of $2 \times 2$ sub-matrices from random orthonormal matrices; Right: empirical pdf for the scores for random $\mathbf{U}_{ \{i,j\} }$ if the entries are standard Gaussian random variables (without the normalization factor $d^{-1/2}$).}
		\label{fig:remark4}
\end{figure}
\begin{figure}[t]
	\centering
	\includegraphics[trim = 45 1 5 5, clip, width=0.85\textwidth]{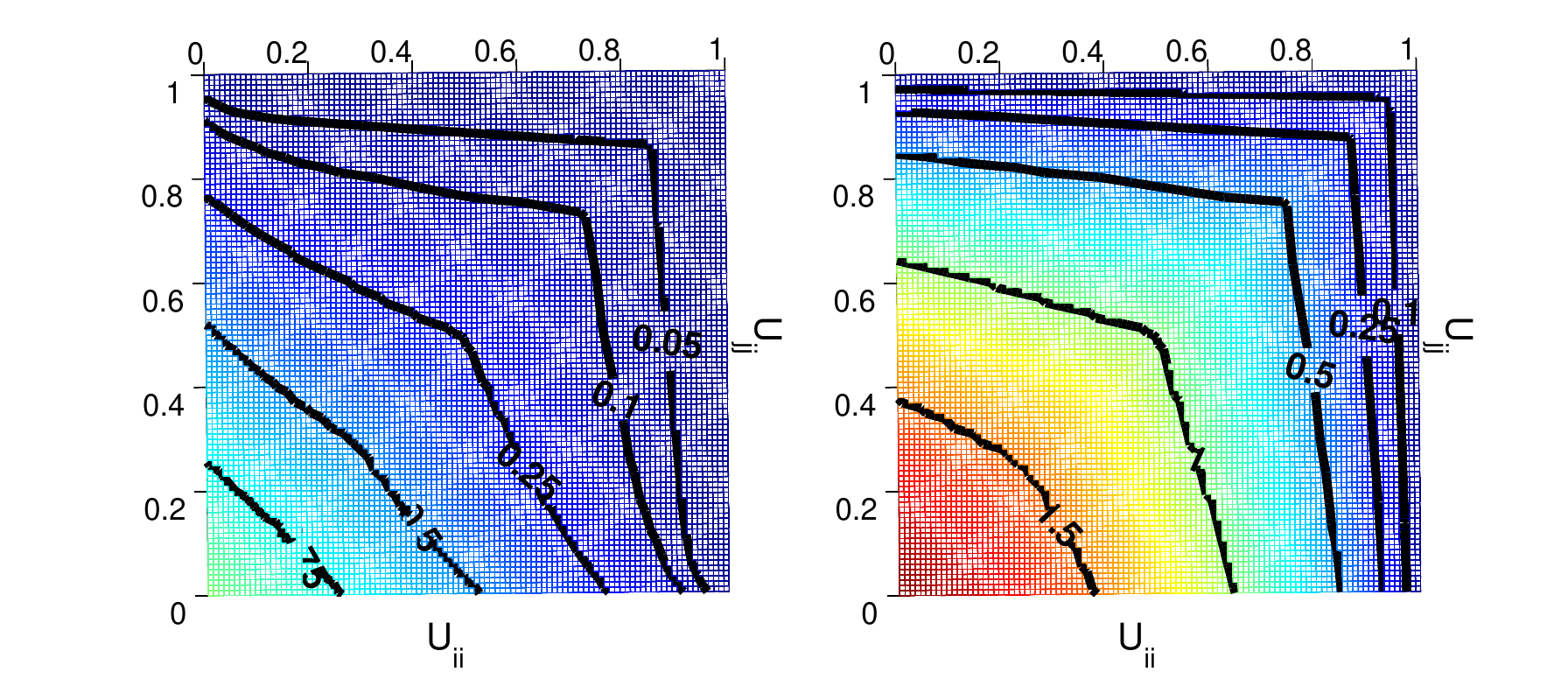}
	\caption{Empirical mean (left) and maximum (right) values of $C_{ij}$ for the toy matrix $\mathbf{U}_{ \{i,j\} }$ where the diagonal elements $U_{ii}$ and $U_{jj}$ (both in $[0,1]$) are fixed and the off-diagonals entries denoted
		%$\mathbf{U}_{ \{i,j\} } = \begin{bmatrix} U_{ii} & z_2 \\ z_1 & U_{jj} \end{bmatrix}$
		$z_1$ and $z_2$ are uniform random variables in the interval $\left[-\sqrt{1-\beta}, \sqrt{1-\beta} \right]$ where $\beta = \max(U_{ii}^2, U_{jj}^2)$ since the columns and rows of $\mathbf{U}$ are $\ell_2$ normalized. Low values (close to 0) are coded as dark blue while high values (close to 2) are coded as dark red.}
	\label{fig:overallaverage}
\end{figure}
\noindent \textbf{More on Remark 4.} Understanding the properties of these $C_{ij}$ is of crucial importance for our approach. First, notice the effect one generalized Givens transformation has: we have $\| \mathbf{U} - \mathbf{G}_{i j} \|_F^2 = \| \mathbf{U}\mathbf{G}_{i j}^T - \mathbf{I}  \|_F^2 = 2d - 2\text{tr}(\mathbf{U}\mathbf{G}_{i j}^T)$, which together with \eqref{eq:givensapproxforortho} leads to
\begin{equation}
\text{tr}(\mathbf{U}\mathbf{G}_{i j}^T) = \text{tr}(\mathbf{U}) + C_{i j} \leq d.
\label{eq:towardsidentity}
\end{equation}
In this sense, the $\mathbf{G}_{i j}$ ``pushes'' $\mathbf{U} \mathbf{G}_{ij}^T$ towards the identity matrix by ``contributing'' $C_{i j}$ to the diagonal of $\mathbf{U}$, i.e., we are estimating the inverse of $\mathbf{U}$ which in this case is just the transpose.\\
Notice that $0 \leq C_{ij} \leq 4$. The minimum is achieved for symmetric positive semidefinite matrices (because in this case the eigenvalues and singular values are the same and therefore the nuclear norm equals the trace) and the maximum for $\mathbf{U}_{\{i, j \} } = -\mathbf{I}_{2 \times 2}$. This immediately leads to a local optimality condition for our approach: there is no $\mathbf{G}_{ij}$ to improve the approximation if all $\mathbf{U}_{ \{ i,j \} }$ are symmetric positive definite. Assume now we are given a random orthogonal matrix. Because the singular values $s_{1,2}$ depend on the entries of $\mathbf{U}_{\{i,j\}  }$ which we model as Gaussian random variables with zero mean and standard deviation $d^{-1/2}$ \cite{Stewart2019}, we have by direct calculation that $\mathbb{E}[s_1] \approx 1.7724 d^{-1/2}$ and $\mathbb{E}[s_2] \approx 0.5190 d^{-1/2}$ which leads to $\mathbb{E}[\| \mathbf{U}_{ \{i,j\} } \|_*] \approx 2.2914 d^{-1/2}$. The trace is the sum of two absolute value Gaussian random variables and therefore $\mathbb{E}[\text{tr}(\mathbf{U}_{ \{i,j\} })] = 2\sqrt{2 (\pi d)^{-1}}$ which leads to
\begin{equation}
\mathbb{E}[C_{ij}] \approx 0.6956 d^{-1/2}.\hfill \square
\label{eq:expectedvalue}
\end{equation}

\noindent \textbf{More on Remark 5.} To see how the scores $C_{ij}$ depend on the diagonal entries of our toy model, by direct calculation with the truncated standard Gaussian random variables we have that
\begin{equation}
\mathbb{E}[C_{ij}(u)] = \frac{1}{2 \pi \text{erf}^2\left( \sqrt{\frac{1-u^2}{2}} \right)}  \iint_{-\sqrt{1-u^2}}^{\sqrt{1-u^2}} \exp \left( -\frac{z_1^2+z_2^2}{2} \right)C_{ij}(u) \ dz_1 dz_2 \propto (1-u)^2.
\end{equation}
We show in Figure \ref{fig:remark4} the empirical results (mean and standard deviation of $C_{ij}$) on the toy matrix $\mathbf{U}_{ \{i,j\}} = \begin{bmatrix}
u & z_2 \\ z_1 & u
\end{bmatrix}$ for $0 < u < 1$. The empirical mean follows the approximation in Remark 5 (and is tight for $u \leq 0.7$) while we notice that the variance is high for almost the whole interval. In Figure \ref{fig:overallaverage} we show the average (left) and maximum (right) costs $C_{ij}$ achieved for another toy model where the diagonal elements are distinct $\mathbf{U}_{\{  i,j \} } = \begin{bmatrix}
U_{ii} & z_2 \\ z_1 & U_{jj}
\end{bmatrix}$ for $0 \leq U_{ii}, U_{jj} \leq 1$.\\
Finally, notice that when $\mathbf{U}_{\{  i,j \} } = \begin{bmatrix}
u & U_{ij} \\ U_{ij} & u
\end{bmatrix}$ we have $C_{ij} = 2(U_{ij} - u)$ if $u \leq U_{ij}$ and zero otherwise, and when $\mathbf{U}_{\{  i,j \} } = \begin{bmatrix}
u & -U_{ij} \\ U_{ij} & u
\end{bmatrix}$ we have $C_{ij} = 2\sqrt{u^2 + U_{ij}^2} - 2u$ indicating that skew symmetric sub-matrices have higher $C_{ij}$ than symmetric ones, in general.$\hfill \square$

\noindent \textbf{Proof of Remark 7.} We use the fact that columns of $\mathbf{U}$ and $\mathbf{\bar{U}}$ have unit norm, and the fact that the Frobenius norm is entrywise:
\begin{equation}
\| (\mathbf{U} - \mathbf{\bar{U}})\mathbf{\Sigma}_d  \|_F^2 = \! \! \sum_{i=1}^d \sigma_i \| \mathbf{u}_i - \mathbf{\bar{u}}_i \|_2^2 = \! \! \sum_{i=1}^d \sigma_i (\|\mathbf{u}_i\|_2^2 + \|\mathbf{\bar{u}}_i\|_2^2 - 2\mathbf{u}_i^T\mathbf{\bar{u}}_i )
= \! \! 2 \! \sum_{i=1}^d \sigma_i (1 - \cos (\theta_i)).\hfill \square
\end{equation}

\noindent \textbf{Proof of Remark 8.} First, note that given any proper norm which is invariant to orthonormal transformations, we have that $\| \mathbf{U} - \mathbf{\bar{U}} \|_2 = \| \mathbf{I} - \mathbf{U}^T\mathbf{\bar{U}} \|_2$ and we denote the error matrix $\mathbf{E} = \mathbf{I} - \mathbf{U}^T\mathbf{\bar{U}}$. Now, start from the error matrix $\mathbf{E} = \mathbf{I} - \mathbf{U}^T \mathbf{\bar{U}}$ and use the fact that $\mathbf{U}^T \mathbf{\bar{U}}$ is orthonormal to notice by straightforward calculation that:
\begin{equation}
\mathbf{E}^T \mathbf{E} = (\mathbf{I} - \mathbf{\bar{U}}^T \mathbf{U})(\mathbf{I} - \mathbf{U}^T \mathbf{\bar{U}}) = \mathbf{I} - \mathbf{U}^T \mathbf{\bar{U}} - \mathbf{\bar{U}}^T \mathbf{U} + \mathbf{I} = \mathbf{E} + \mathbf{E}^T.
\end{equation}
Denote now a complex-valued eigenvalue of $\mathbf{E}$ by $z_k = a_k + ib_k$, then because of the equality above we have that $z_{k}^* z_k = 2 \Re(z_k)$ which in turn can be written as $a_k^2 + b_k^2 = 2a_k$ or $(a_k-1)^2 + b_k^2 = 1$. This is the equation of a circle of radius one centered at $(1,0)$ in the complex plane. Another way to view this result is to observe that $\mathbf{U}^T \mathbf{\bar{U}}$ is orthonormal an therefore its spectrum is on the unit circle and then $\mathbf{I}$ shifts the whole spectrum to the right by one unit.\\
Moreover, we also have that $\mathbf{EE}^T = \mathbf{E} + \mathbf{E}^T$ and therefore $\mathbf{E}^T \mathbf{E} = \mathbf{EE}^T$ which means that $\mathbf{E}$ is a normal matrix. As such, the singular values of $\mathbf{E}$ are the absolute values of its eigenvalues and therefore $\| \mathbf{E} \|_2 \leq 2$.$\hfill \square$
\newline

\noindent \textbf{Proof of Theorem 5.} The proof is based on the Gershgorin circle theorem (detailed in Chapter 7.2 of \cite{Golub1996}) applied to the error matrix $\mathbf{E} = \mathbf{I} - \mathbf{U}^T \mathbf{\bar{U}}$. If we denote $\mathbf{Q} = \mathbf{U}^T \mathbf{\bar{U}}$, we then have for the $i^\text{th}$ eigenvalue of $\mathbf{E}$ that
\begin{equation}
| \lambda_i - 1 + \mathbf{u}_i^T \mathbf{\bar{u}}_i  | \leq \sum_{j\neq i} |Q_{ij}| \leq \sqrt{(d-1)(1-(\mathbf{u}_i^T \mathbf{\bar{u}}_i)^2)}.
\label{eq:lambdanow}
\end{equation}
The last inequality on the right hand side comes from the $\ell_1-\ell_2$ inequality $\| \mathbf{x} \|_1 \leq \sqrt{d-1}\| \mathbf{x} \|_2$ applied to the vectors of size $d-1$ (the rows of $\mathbf{Q}$ except for their diagonal elements and whose $\ell_2$ norm is $\sqrt{1-Q_{ii}^2}$). Returning to \eqref{eq:lambdanow}, for the term on the left hand side we have by the reverse triangle inequality ($|a-b| \geq | |a| - |b| |$) that
\begin{equation}
| |\lambda_i| - |1 - \mathbf{u}_i^T \mathbf{\bar{u}}_i| | \leq | \lambda_i - (1 - \mathbf{u}_i^T \mathbf{\bar{u}}_i)  |.
\label{eq:reversetriangle}
\end{equation}
Finally combining \eqref{eq:lambdanow} and \eqref{eq:reversetriangle} we have that
\begin{equation}
0 \leq |\lambda_i| \leq 1- \mathbf{u}_i^T \mathbf{\bar{u}}_i + \sqrt{(d-1)(1-(\mathbf{u}_i^T \mathbf{\bar{u}}_i)^2)}.
\end{equation}
%|x - |1-c|| <= sqrt((d-1)*(1-c^2)), c<=1, d>=2, x>=0 in wolfram
This inequality holds for $0 \leq \mathbf{u}_i^T \mathbf{\bar{u}}_i \leq 1$ and $d \geq 2$. The expression on the right-hand side is monotonic decreasing for $\mathbf{u}_i^T \mathbf{\bar{u}}_i \geq 0$. As such, the largest upper bound happens for $i = \underset{i}{\arg \min}\ \mathbf{u}_i^T \mathbf{\bar{u}}_i$. And because, as shown in Remark 8, the error matrix is a normal matrix we get a bound on the operator norm (the singular values of a normal matrix are the absolute values of its eigenvalues). $\hfill \blacksquare$

\end{document}